\documentclass[11pt]{article}

\usepackage{amssymb,diagrams}
\usepackage{epsfig,amsmath,latexsym}
\usepackage{amsfonts}

\newtheorem{theorem}{Theorem}[section]

\newtheorem{defi}{Definition}[section]

\newtheorem{prop}[theorem]{Proposition}

\def\slfrac#1#2{\hbox{\kern.1em %
 \raise.5ex\hbox{\the\scriptfont0 #1}\kern-.11em %
 /\kern-.15em\lower.25ex\hbox{\the\scriptfont0 #2}}}

\newarrow{Into}C--->
\newarrow{Line}-----

\newcommand{\pf}{\noindent{\bf Proof.~}}
\newcommand{\beq}{\begin{eqnarray}}
\newcommand{\eeq}{\end{eqnarray}}
\newcommand{\beql}[1]{\begin{eqnarray}\label{#1}}
\newcommand{\beqs}{\begin{eqnarray*}}
\newcommand{\eeqs}{\end{eqnarray*}}
\newcommand{\eqn}[1]{(\ref{#1})}

\newcommand{\CC}{{\mathbb C}}
\newcommand{\DD}{{\mathbb D}}
\newcommand{\HH}{{\mathbb H}}

\newcommand{\PP}{{\mathbb P}}

\newcommand{\RR}{{\mathbb R}}

\newcommand{\cc}{{\mathbb C}}
\newcommand{\rr}{{\mathbb R}}
\newcommand{\zz}{{\mathbb Z}}

\newcommand{\dd}{{\mathbb D}}

\newcommand{\fc}{{\mathfrak c}}
\newcommand{\fd}{{\mathfrak d}}

\newcommand{\fg}{{\mathfrak g}}

\newcommand{\fj}{{\mathfrak j}}

\newcommand{\fp}{{\mathfrak p}}

\newcommand{\fs}{{\mathfrak s}}
\newcommand{\ft}{{\mathfrak t}}

\newcommand{\bb}{{\mathbf b}}

\newcommand{\br}{{\mathbf r}}
\newcommand{\be}{{\mathbf e}}

\newcommand{\bh}{{\mathbf h}}
\newcommand{\bm}{{\mathbf m}}
\newcommand{\bo}{{\mathbf 1}}
\newcommand{\bs}{{\mathbf s}}

\newcommand{\bw}{{\mathbf w}}
\newcommand{\bx}{{\mathbf x}}

\newcommand{\bz}{{\mathbf z}}

\newcommand{\bA}{{\mathbf A}}

\newcommand{\bD}{{\mathbf D}}
\newcommand{\bI}{{\mathbf I}}
\newcommand{\bJ}{{\mathbf J}}

\newcommand{\bM}{{\mathbf M}}

\newcommand{\bP}{{\mathbf P}}
\newcommand{\bQ}{{\mathbf Q}}
\newcommand{\bS}{{\mathbf S}}
\newcommand{\bT}{{\mathbf T}}
\newcommand{\bU}{{\mathbf U}}
\newcommand{\bV}{{\mathbf V}}
\newcommand{\bW}{{\mathbf W}}

\newcommand{\bY}{{\mathbf Y}}
\newcommand{\bZ}{{\mathbf Z}}

\newcommand{\sA}{{\mathcal A}}

\newcommand{\sD}{{\mathcal D}}

\newcommand{\sG}{{\mathcal G}}

\newcommand{\sL}{{\mathcal L}}
\newcommand{\sM}{{\mathcal M}}
\newcommand{\sP}{{\mathcal P}}
\newcommand{\sS}{{\mathcal S}}

\newcommand{\sX}{{\mathcal X}}
\newcommand{\PD}{\sP_{\sD}}

\newcommand{\bsq}{\vrule height .9ex width .8ex depth -.1ex}

\makeatletter
\def\@sect#1#2#3#4#5#6[#7]#8{\ifnum #2>\c@secnumdepth
     \def\@svsec{}\else
     \refstepcounter{#1}\edef\@svsec{\csname the#1\endcsname.\hskip .75em }\fi
     \@tempskipa #5\relax
      \ifdim \@tempskipa>\z@
        \begingroup #6\relax
          \@hangfrom{\hskip #3\relax\@svsec}{\interlinepenalty \@M #8\par}%
        \endgroup
       \csname #1mark\endcsname{#7}\addcontentsline
         {toc}{#1}{\ifnum #2>\c@secnumdepth \else
                      \protect\numberline{\csname the#1\endcsname}\fi
                    #7}\else
        \def\@svsechd{#6\hskip #3\@svsec #8\csname #1mark\endcsname
                      {#7}\addcontentsline
                           {toc}{#1}{\ifnum #2>\c@secnumdepth \else
                             \protect\numberline{\csname the#1\endcsname}\fi
                       #7}}\fi
     \@xsect{#5}}
\def\@begintheorem#1#2{\it \trivlist \item[\hskip \labelsep{\bf #1\ #2.}]}

\def\plain{plain}\ifx\fmtname\plain\csname fi\endcsname
     
     \input docstrip
     \preamble

     Do not distribute the stripped version of this file.
     The checksum in the header refers to the documented version.

     \endpreamble
     \generateFile{here.sty}{t}{\from{here.doc}{}}
     \endinput
\fi
\ifcat a\noexpand @\let\next\relax\else\def\next{%
    \documentstyle[here,doc]{article}\MakePercentIgnore}\fi\next
\ifx\@Hxfloat\@Hundef\else\expandafter\endinput\fi
\let\@Hxfloat\@xfloat
\def\@xfloat#1[{\@ifnextchar{H}{\@HHfloat{#1}[}{\@Hxfloat{#1}[}}
\def\@HHfloat#1[H]{%
\expandafter\let\csname end#1\endcsname\end@Hfloat
\vskip\intextsep\vbox\bgroup\def\@captype{#1}\parindent\z@
\ignorespaces}
\def\end@Hfloat{\egroup\vskip \intextsep}

\makeatother

\catcode`\@=11

\catcode`@=12

\begin{document}


\begin{center}
{\Large {\bf Apollonian Circle Packings: Geometry and Group Theory \\
I. The Apollonian Group}}\\

\vspace{1.5\baselineskip}
{\em Ronald L. Graham} \\
Department of Computer Science and Engineering \\
 University of California at San Diego, 
La Jolla, CA 92093-0114 \\
\vspace*{1.5\baselineskip}

{\em Jeffrey C. Lagarias} \\
Department of Mathematics \\
University of Michigan, 
Ann Arbor, MI 48109--1109 \\
\vspace*{1.5\baselineskip}

{\em Colin L. Mallows} \\
Avaya Labs, Basking Ridge, NJ 07920 \\
\vspace*{1.5\baselineskip}

{\em Allan R. Wilks} \\
AT\&T Labs, Florham Park, NJ 07932-0971 \\
\vspace*{1.5\baselineskip}

{\em Catherine H. Yan}
\footnote{
Partially supported by NSF grants DMS-0070574, DMS-0245526 and a Sloan
Fellowship. This author is also affiliated with Dalian University of
Technology, China.}\\
Department of Mathematics \\
Texas A\&M University, College Station, TX 77843-3368\\
\vspace*{1.5\baselineskip}
\vspace*{1.5\baselineskip}
(March 10, 2005  version) \\
\vspace*{1.5\baselineskip}
{\bf ABSTRACT}
\end{center}
Apollonian circle packings arise by repeatedly filling the interstices
between four mutually tangent circles with further tangent circles.
We observe that there exist 
Apollonian  packings which have strong integrality
properties, in which  all  circles in the packing have
 integer curvatures and
rational centers such that (curvature)$\times$(center) is an
integer vector. This series of papers explain such properties.

A {\em Descartes configuration} is a set of four mutually tangent
circles with disjoint interiors. 
 An Apollonian circle packing can be described in terms of
the Descartes configuration it contains. We describe the 
 space of all  ordered, oriented Descartes configurations
using  a coordinate system $\sM_\DD$
consisting of  those $4 \times 4$
real matrices $\bW$ with $\bW^T \bQ_{D} \bW = \bQ_{W}$ where 
$\bQ_D$ is the matrix of the Descartes quadratic form
$Q_D= x_1^2 + x_2^2+ x_3^2 + x_4^2 - 
\frac{1}{2}(x_1 +x_2 +x_3 + x_4)^2$ 
and  $\bQ_W$ of the quadratic form
$Q_W = -8x_1x_2 + 2x_3^2 + 2x_4^2$. 
On the parameter space
$\sM_\DD$ the  group
$Aut(Q_D)$ acts on the left, and $Aut(Q_W)$ acts on the right,
giving two different ``geometric'' actions.
Both these  groups are  isomorphic to the 
Lorentz group $O(3, 1)$. 
The right action of $Aut(Q_W)$ 
(essentially) corresponds to Mobius transformations acting on the underlying
Euclidean space $\rr^2$ while the left action of $Aut(Q_D)$ is defined
only on the parameter space. 

We observe that
the Descartes configurations in each Apollonian packing 
form  an orbit of a single Descartes
configuration under a certain finitely generated discrete
subgroup of  $Aut(Q_D)$, 
which we call the {\em Apollonian group}.
This group consists of 
$4 \times 4$ integer matrices, and its integrality properties 
lead to the integrality
properties observed in some Apollonian circle packings.

We introduce two more related finitely generated groups
in   $Aut(Q_D)$, the {\em  dual Apollonian group} produced 
from the Apollonian
group by a ``duality'' conjugation, and the {\em super-Apollonian group} 
which is the group generated by the Apollonian and
dual Apollonian groups together. 
These groups also
consist of integer $4 \times 4$ matrices. We show
these groups are hyperbolic Coxeter groups.

\vspace*{1.5\baselineskip}
\newpage
\setcounter{page}{1}

%
%
%
\setlength{\baselineskip}{1.0\baselineskip}

\section{Introduction}
\setcounter{equation}{0}
An Apollonian circle packing is a packing of circles 
arising  by repeatedly filling the interstices
between four mutually tangent circles with further tangent circles.
We call an initial arrangement of four mutually tangent circles with
distinct tangents (necessarily $6$ of them) a
{\em Descartes configuration}.

%
%
%

\begin{figure}[htbp]
\centerline{\epsfxsize=4.0in \epsfbox{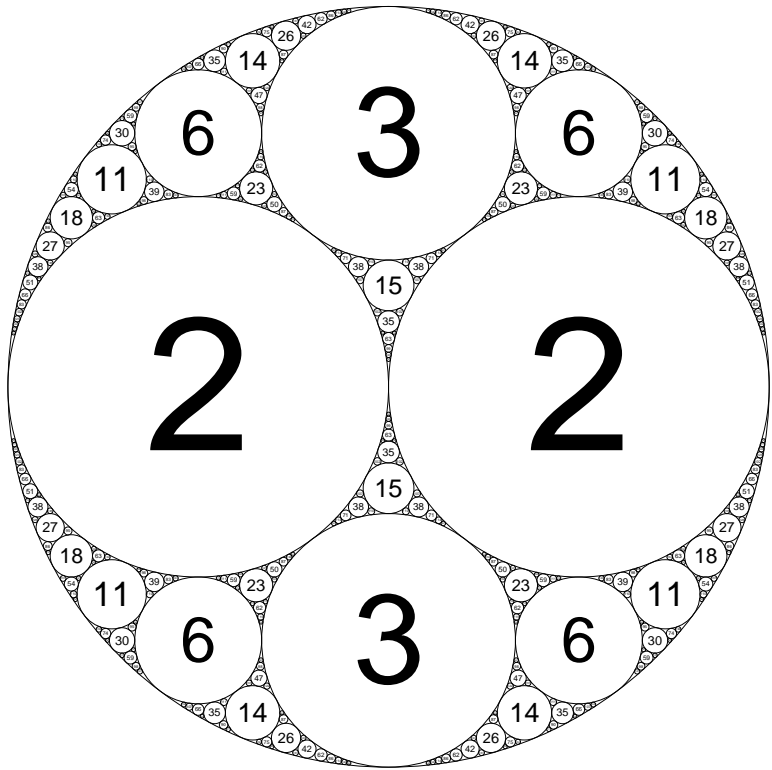}}
\caption{An Apollonian packing, labelled with circle curvatures}\label{std}
\end{figure}

Starting from any Descartes configuration,
we  can recursively construct an  infinite circle packing of the 
Euclidean plane, 
 in which new circles are added 
which are tangent to three of the circles that have already been
placed and have interiors disjoint from
any of them. The infinite packing obtained in the limit of adding
all possible such circles is called an
{\em Apollonian packing}. The new circles added
at each stage can be
obtained using M\"{o}bius transformations of 
Descartes configurations in the partial packing. 

An Apollonian packing is pictured in Figure \ref{std}, in
which each circle is labelled by its {\em curvature,} which is
the inverse of its radius.  The initial
configuration consists of two circles of radius $\frac{1}{2}$
inscribed in a circle of radius $1$, the latter being assigned
negative curvature $-1$ by a convention given in \S3.
This particular Apollonian
packing has the special property that all circles in the packing
have integer curvatures. We call a packing with this
property {\em integral.} More remarkably,
if one regards the circle centers 
as complex numbers, then one can place the initial circles so
that every circle in the packing has ``curvature$\times$center''
a Gaussian integer (element of the ring $\zz[i].)$ This occurs 
for example when the
center of the outer circle is placed at $\bz= 0$ and the
two circles of radius $\frac{1}{2}$ have centers at $\bz= -\frac{1}{2}$
and at $\bz= \frac{1}{2}$. We call a packing
with this extra property {\em strongly integral}. 

The object of this paper is to give a  geometric explanation of
the origin of these integrality properties involving both
the curvatures and the circle centers.
This is based on five facts:

(1) An Apollonian packing can be described in terms of
the Descartes configurations it contains.

(2) There is a coordinate representation of the set $\sM_{\DD}$
of all ordered, oriented Descartes configurations 
as a six-dimensional  real-algebraic variety.
It consists of  the set of  $4 \times 4$ real matrices $\bW$ satisfying a
system of quadratic equations, $\bW^T \bQ_{D} \bW = \bQ_{W}$,
which state that $\bW$ conjugates 
the Descartes quadratic form $Q_{D}= x_1^2 + x_2^2+ x_3^2 + x_4^2 - 
\frac{1}{2}(x_1 +x_2 +x_3 + x_4)^2$
to a quadratic  form 
$Q_{W}= -8x_1x_2 + 2x_3^2 + 2x_4^2$ that
we call the Wilker quadratic form, after Wilker \cite{Wi81}. 
This coordinate system we call ``augmented curvature-center
coordinates'', 
as it encodes the curvatures and centers
of the circles in the configuration. We term
 $\sM_{\DD}$ with these coordinates the {\em parameter space} of 
ordered, oriented Descartes configurations.

(3) The variety $\sM_{\DD}$ is a principal
homogeneous space for the Lorentz group 
$O(3,1)$, under both a left and a right action of this
group, realized as $Aut(Q_{D})$ and $Aut(Q_{W}),$
respectively. The right action corresponds to M\"{o}bius
transformations acting on the plane, while  the left
action acts only on the Descartes configuration space.  

(4) There is a discrete subgroup of the left action
by $Aut(Q_{D})$,  the {\em Apollonian group} $\sA$,
having the property that the (unordered, unoriented) 
Descartes configurations
in any Apollonian packing are described by a single
orbit of this group. If the Descartes configurations
are regarded as ordered and oriented, then 
exactly $48$ such orbits correspond to each packing,
each one containing a copy of each unordered,
unoriented Descartes configuration in the packing.

(5) The Apollonian group $\sA$ consists of  integer
$4 \times 4$ matrices.

The last property (5) explains the existence of Apollonian
packings having integrality properties of both curvatures
and centers. If an initial Descartes configuration has
integral augmented curvature-center coordinates, then the 
same property persists for all Descartes configurations in the
packing, hence for all circles in the packing. The packing
pictured in Figure~\ref{std} has this property.

The observation that there are Apollonian packings having
all curvatures integral is an old one, and was noted by
F. Soddy \cite{Sod36}, \cite{Sod37}, who also gave  an
extension to three dimensions, the ``bowl of integers''.
The existence of an integral matrix group  explaining such  curvatures
was studied in S\"{o}derberg \cite{So92} in 1992. 
There are many other mentions in the literature. What we 
add in this direction  is the observation that the circle centers
can simultaneously have integrality properties. 

A main contribution of these papers is the viewpoint 
that these properties should be understood in
terms of group actions on the
parameter space $\sM_{\DD}$ of all 
ordered, oriented Descartes configurations.
Besides introducing 
coordinate systems for this parameter space
(given first in \cite{LMW02} by three of the authors), 
we  study the relevant discrete group actions in detail.
In particular we introduce a larger group of integer
matrices acting on the left, the {\em super-Apollonian group}
$\sA^{S}$, whose generators have natural geometric interpretations
in terms of their action on Descartes configurations. 
We prove that this group is a hyperbolic Coxeter group.

The detailed contents of part I are summarized in
the next section. Below we briefly
indicate the contents of parts II and III.

In part II we study the integrality properties
of the Apollonian packing in more detail, 
and their relation to the super-Apollonian group. 
We show that every integral Apollonian packing can
be transformed to a strongly integral one by a Euclidean
motion. We introduce super-packings as orbits of
the super-Apollonian group, starting
from a given Descartes configuration, and geometric
super-packings which consist
of  the set of circles in all these
Descartes configurations. We show that there 
are exactly 8 different primitive
\footnote{An integral packing is {\em primitive} if
the greatest common divisor of all curvatures of 
circles in the packing is $1$.}
strongly integral geometric super-packings,
and that each of these contains a copy of every
primitive integral Apollonian packing. 
We  characterize
the set of all (primitive) strongly integral Descartes configurations
as a collection of $384$ orbits of the super-Apollonian group.

In part III we consider to what extent the results
proved in parts I and II extend to higher dimensions.  
In all dimensions $n \ge 3$ there are analogues
of the parameter space of Descartes configurations, of
left and right group actions by $O(n+1, 1)$, and of
the Apollonian, dual Apollonian  and super-Apollonian groups. 
However in dimensions $n \ge 4$ the Apollonian group
action does not correspond to a sphere packing; the
spheres will overlap. Nevertheless one can still study the
orbits of these group actions on the parameter space.
The Apollonian group in higher dimensions has rational
entries, rather than integer entries. We show that 
configurations having all curvatures rational exist in
all dimensions, and having  curvature$\times$centers
rational exist in certain dimensions only.

The general framework of these papers was developed
by the second author (JCL), who also did much of
the writing. This paper is an extensively  revised version of a preprint
written in 2000, which adds some new results in \S4 and S5. 

\paragraph {Acknowledgments.}
The authors thank Ludwig Balke for suggesting that the 
two group actions in \S3 should correspond to a 
left action and a right action. The authors are  grateful for
helpful comments from Andrew Odlyzko, Eric Rains, 
Jim Reeds and Neil Sloane.  The authors thank the reviewer
for incisive comments leading to reorganization of the paper.

%
%
%

\section{Summary of  Results}
\setcounter{equation}{0}

There are at least three  ways to 
describe the Apollonian packing $\sP$ containing a given
Descartes configuration $\sD$. \\

(G1) [Geometric] An Apollonian  packing $\sP_{\sD}$ 
is a set of circles in 
$\hat{\cc} = \RR^2 \cup \{ \infty\}$, which comprise four orbits
under the action of a discrete group $G_{\sA}(\sD)$ of
M\"{o}bius transformations
inside  the conformal group M\"{o}b(2). The discrete
group $G_{\sA}(\sD)$ depends on $\sD$. \\

(G2) [Algebraic] The collection of all (ordered, oriented) Descartes
configurations in the packing $\sP_{\sD}$ form $48$ orbits
of a discrete group $\sA$, the Apollonian group, contained
in the group $Aut(\sM_{\DD}) \equiv Aut(Q_{D})$ of
left-automorphisms of the parameter space $\sM_{\DD}$ of Descartes
configurations. The discrete group $\sA$ is independent of 
the configuration $\sD$. \\

(G3) [Holographic]   
The open disks comprising the interiors of the
circles in the packing are the complement
$\hat{\cc} \backslash \Lambda_{\sD}$ of the limit
set $\Lambda_{\sD}$ of a certain Schottky group $\sS_{\sD}$
acting on hyperbolic $3$-space $\HH^3$, with  $\hat{\cc}$
identified with  its ideal boundary. 
The Schottky group  $\sS_{\sD}$ depends on $\sD$. \\

In this paper we will mainly consider viewpoints
(G1) and (G2). Viewpoint (G3) is described in 
Mumford, Series and Wright \cite[Chapter 7]{MSW02},
and we treat it in  Appendix B. It is termed
``holographic'' because it views the 
limit set of the packing as 
the boundary of a higher-dimensional object, 
which in principle  gives information about it. In
Appendix B we point out 
a connection of viewpoint (G3) to viewpoint (G1).
The term ``algebraic'' for (G2) refers to the group action
being on a real-algebraic variety  of Descartes 
configurations.

The main emphasis of this series of papers 
is to study Apollonian packings in terms of the 
Descartes configurations they contain.
That is, we  study the packing
as a collection of points
 inside the parameter space of all 
ordered, oriented Descartes configurations. 
There are  two different group actions on this parameter space,
which are a right action associated to viewpoint (G1)
and a left action corresponding to viewpoint (G2), 
with the groups both isomorphic to
the Lorentz group $O(3,1)$, a real Lie group, as explained
in \S3. In particular  group actions are independent and 
mutually commute.
The discrete group $G_{\sA}(\sD)$ above is contained in 
the right action and the discrete group $\sA$, the Apollonian
group,  is contained in the  left action.
Thus viewpoints (G1) and (G2) are complementary and
coexist simultaneously on the space of all Descartes configurations.

As stated in the introduction,
the integer structures in the curvatures
and centers of some Apollonian packings can be explained
in terms of the viewpoint (G2),
using the discrete group $\sA$,
the Apollonian group.  This viewpoint can be traced  back to
the ``inversive crystal'' in Wilker~\cite[\S 14]{Wi81}.

In \S3 we coordinatize 
the space  $\sM_{\DD}$ of
(ordered, oriented) Descartes configurations as in \cite{LMW02}
and describe two group actions on this space.
In \S3.1 we give two coordinate systems.
The first labels such a Descartes configuration with a 
$4 \times 3$ matrix $\bM_{\sD}$, 
called {\em curvature-center coordinates}, 
and 
the second with a $4 \times 4$ matrix $\bW_{\sD}$,
called {\em augmented curvature-center coordinates}. 
These coordinates
are characterized  by 
quadratic relations which generalize the Descartes circle
theorem. Theorem~\ref{th31} exactly characterizes
these relations for curvature-center coordinates, 
strengthening Theorem 3.2 in \cite{LMW02}
by formulating and proving a converse.
For the augmented curvature-center coordinates
these quadratic relations (proved in  \cite{LMW02})
take the form
$$
\bW_{\sD}^{T} \bQ_{D} \bW_{\sD} = \bQ_{W},
$$
which gives  a conjugacy  of the Descartes quadratic form $Q_{D}$
to the Wilker quadratic form $Q_{W}$.
These forms are indefinite of signature $(3,1)$ and their
(real) automorphism groups $Aut(Q_{D})$ and $Aut(Q_W)$
are isomorphic to the Lorentz group $O(3,1)$.

In \S3.2 we describe  a linear  left action (``Lorentz action'')
by $Aut(Q_{D})$ and a right action (``M\"{o}bius action'') 
by $Aut(Q_W)$ on the space  $\sM_{\DD}$.
Theorem~\ref{th33} describes these actions. It shows
that the space $\sM_{\DD}$ is a principal homogeoneous
space for the group $O(3,1)$ under both the left action
and the right action. 
The M\"{o}bius action
is treated in more detail in Appendix A.
In \S3.3  we describe some integral elements of 
the Lorentz action $Aut(Q_{D})$ which have geometric
interpretations as simple transformations  of a Descartes
configuration. These elements  are used in defining the 
Apollonian group, dual Apollonian group and super-Apollonian
group given below.

In \S4 we describe Apollonian packings and the 
Apollonian group. Theorem~\ref{th41} establishes
the basic fact that the interiors of all circles
in an Apollonian packing are disjoint.
The much-studied {\em residual set} $\Lambda(\sP)$ of an
Apollonian packing is the complement of the
interiors of all circles; it is a set of 
measure zero.  For later use, 
Theorem~\ref{th42} gives several properties  of  the
residual set $\Lambda(\sP)$.
Theorem~\ref{th41}  shows that the 
(ordered, oriented) Descartes configurations
in a packing form 48 orbits of this group.

In \S5 we define the dual Apollonian group $\sA^{\perp}$
and call its orbits {\em dual Apollonian packings}.
Theorem \ref{th43} shows that the  set of all
tangency points of circles in a dual Apollonian packing
has closure the limit set of another Apollonian packing,
that generated by the dual Descartes configuration of
any Descartes configuration generating the packing. 

In \S6 we define the super-Apollonian group $\sA^{S}$ to
be the group generated by $\sA$ and $\sA^{\perp}$ combined,
and call its orbits Apollonian super-packings.
It is a discrete group of integer matrices contained in $Aut(Q_D)$.
Theorem~\ref{Sth47} in \S6.1
gives  a complete  presentation for $\sA^{\perp}$, 
establishing that it is a hyperbolic Coxeter group.
In \S6.2 we add remarks on  super-packings, which are
studied at length in part II.

In Appendix A we describe the M\"{o}bius group action in
detail. Theorem~\ref{th91} gives an isomorphism of this group to 
$Aut(Q_{W})$. 

In Appendix B we describe  the Schottky group
action in Mumford, Series and Wright \cite{MSW02}.
We indicate some relations to the  
M\"{o}bius group action.

%
%
%

\section{Descartes Configurations and Group Actions}
\setcounter{equation}{0}

%
%

\subsection{Descartes Configurations and Curvature-Center Coordinates}

In 1643, Descartes found a relation between the
radii for four mutually disjoint tangent circles of type (a)
in Figure \ref{fig1} below.
Let $r_1, r_2, r_3, r_4$ be the
radii of the tangent circles. 
Descartes showed~\footnote{Descartes expressed his relation in
a different form,  obtained by clearing denominators in \eqn{001}.}
a result equivalent to 
\beql{001}
\frac{1}{r_1^2}+\frac{1}{r_2^2}+\frac{1}{r_3^2}+\frac{1}{r_4^2}
=\frac{1}{2} \left(\frac{1}{r_1}+
\frac{1}{r_2}+\frac{1}{r_3} +\frac{1}{r_4}\right)^2,
\eeq
which is now called the {\em Descartes circle theorem}.
This result can be extended to apply to all Descartes
configurations, which can be of types (a)-(d) in
the figure below.

%
%

\begin{figure}[htbp]
\centerline{\epsfxsize=3.0in \epsfbox{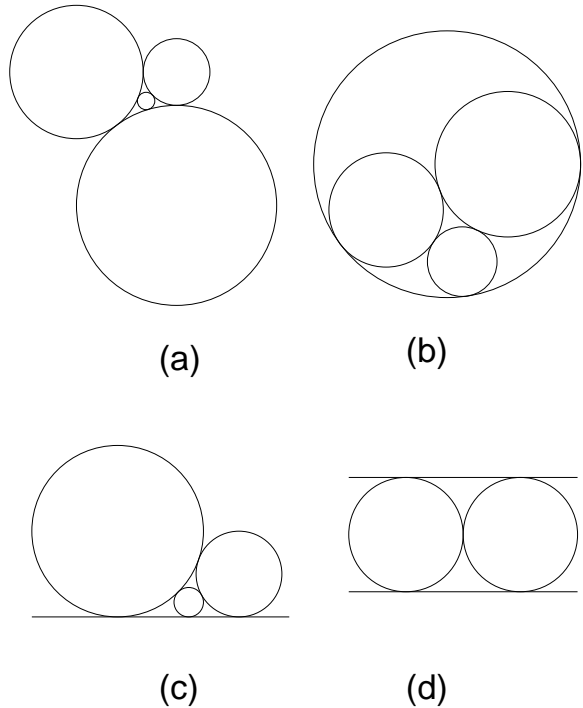}}
\caption{Descartes configurations}\label{fig1}
\end{figure}

In Figure \ref{fig1} (a) is a standard
Descartes  configuration in which the circles have disjoint interiors,
and all curvatures are positive.
We will also allow configurations like Figure~\ref{fig1} (b) 
in which one circle
encloses the other three.
We also allow straight lines to be
regarded as circles of curvature zero, leading to Descartes configurations
like Figure~\ref{fig1} (c) and (d), in which the ``interior'' 
of the ``circle''
defined by a straight line is defined to be a suitable half-plane.

To extend the Descartes Circle Theorem to all Descartes configurations
of types (a)-(d), we must define the curvatures to have
appropriate signs, as follows.  An {\em oriented circle} is a circle
together with an assigned direction of unit normal vector, which can
point inward or outward.  If it has radius $r$ then its {\em oriented
radius} is $r$ for an inward pointing normal and $-r$ for an outward
pointing normal.  Its {\em oriented curvature} (or ``signed
curvature'') is $1/r$ for an inward pointing normal and $-1/r$ for an
outward pointing normal.  By convention, the {\em interior} of an
oriented circle is its interior for an inward pointing normal and its
exterior for an outward pointing normal.  An {\em oriented Descartes
configuration} is a Descartes configuration in which the orientations
of the circles are compatible in the following sense:  either (i) the
interiors of all four oriented circles are disjoint, or (ii) the
interiors are disjoint when all orientations are reversed.  Each
Descartes configuration has exactly two compatible orientations in
this sense, one obtained from the other by reversing all orientations.
The {\em positive orientation} of a Descartes configuration is
the one in which the sum of the signed curvatures is positive, while
the {\em negative orientation} is the one in which the sum of
the curvatures is negative. (One can show 
that the sum of the signed curvatures is always nonzero.)  
Note that positive orientation corresponds to case (i),
and negative orientation to case (ii). With these
definitions, the Descartes circle theorem remains valid for all
oriented Descartes configurations, using oriented curvatures.

The geometry of Descartes configurations in Figure \ref{fig1}
is encoded in the curvature
vector $\bb=(b_1,b_2,b_3,b_4)^T$ where $b_i$ is the oriented curvature 
of the $i$th circle. If $\sum_{j = 1}^{4} b_j > 0$, then one of the
following holds:  (a) all of $b_1, ~b_2, b_3, b_4$ are positive;
(b) three are positive and one is negative; (c) three are positive
and one is zero; or (d) two are positive and equal and the other two
are zero. 

The Descartes circle theorem gives a quadratic equation
for the curvature vector $\bb$, which can be rewritten as
\beql{002}
\bb^T \bQ_{D} \bb = 0,
\eeq
in which
\beql{003}
\bQ_{D}= \bI - \frac{1}{2}{\bf 1}{\bf 1}^T =          
         \frac{1}{2}\left[ \begin{array}{rrrr}
               1 & -1 & -1 & -1 \\
               -1 & 1 & -1 & -1 \\
               -1 & -1 & 1 & -1 \\
               -1 & -1 & -1 & 1 \end{array} \right]   
\eeq
is the {\em Descartes quadratic form.} Here $\bo$ represents the vector
$(1,1,1,1)^T$. 

In \cite{LMW02} three of the authors of this paper
showed there exist matrix extensions of the
Descartes circle theorem which encode information
about both the curvatures and the centers of the
circles in the configuration, as follows.
Given an oriented circle $C$ with center $(x, y)$
and oriented curvature $b$  we define its
{\em curvature-center coordinates}  to be the $1 \times 3$
row vector $\bm(C) := (b, bx, by)$.
For the ``degenerate case'' of  an oriented  straight line $H$ we 
define its curvature-center coordinates as 
\[
\bm(H) :=(0, \bh),
\]
where $\bh=(h_1, h_2)$ is the unit normal vector giving the orientation of 
the straight line. 
%
%
%
%

\begin{theorem}[Extended Descartes Theorem]~\label{th31}
Given an ordered configuration $\sD$ 
of four oriented circles with 
oriented 
curvatures $(b_1,~b_2,~b_3,~b_4)$ and centers
$\{(x_i, y_i)~: 1 \leq i \leq 4\}$,
let $\bM_\sD$ be the $4 \times 3$ matrix
\beql{208}
\bM_\sD :=
        \left[ \begin{array}{ccl}
           b_1 & b_1x_1 & b_1y_1 \\
           b_2 & b_2x_2 & b_2y_2 \\
           b_3 & b_3x_3 & b_3y_3 \\
           b_4 & b_4x_4 & b_4y_4
        \end{array} \right] .
\eeq
We include the ``degenerate cases'' where some circles are
oriented lines. 
If this  configuration is an oriented Descartes configuration, then
$\bM=\bM_\sD$ has a nonzero first column and satisfies  
\beql{209}
        \bM^T \bQ_{D} \bM =
        \left[ \begin{array}{ccc}
           0 & 0 & 0 \\
           0 & 2 & 0 \\
           0 & 0 & 2
        \end{array} \right].
\eeq
Conversely, any real solution  $\bM$ to \eqn{209} 
with a nonzero first column is the 
curvature-center coordinate matrix $\bM_{\sD}$ of a unique 
ordered, oriented Descartes configuration $\sD$.  
\end{theorem}

\paragraph{Remark.} The 
hypothesis of a nonzero first column is
necessary in Theorem~\ref{th31}.
Given  the matrix $\bM$ of a Descartes configuration, the
matrix  $\tilde{\bM}$ obtained by zeroing out its first
column will continue to satisfy  \eqn{209}, and will not come from
a Descartes configuration.

\paragraph{Proof.}
An $n$-dimensional version of this theorem was 
stated as Theorem 3.2 in \cite{LMW02}, and the
``if'' direction of the result was proved there. 
(The converse part of the theorem was not
proved there.) More precisely, in  \cite{LMW02} it was proved
that the curvature-center
coordinate matrix $\bM$
of any ordered, oriented  Descartes configuration
satisfies \eqn{209}.   The  first
column of the matrix $\bM$ of a Descartes
configuration cannot be identically zero because at
least two circles in any Descartes configuration have 
nonzero curvature.

To establish the  converse part of the theorem,
we need Theorem \ref{th32} below, which was
independently proved in \cite{LMW02}.  We postpone the proof 
to the end of \S3.1. 
$~~~\bsq$

%
%

The  curvature-center coordinates $\bm(C)$  uniquely determine
a circle in general position, but they fail to be unique
for the degenerate case of a line, where the information
$(0, \bh)$ determines only a  direction normal to the line
and not its position.

In  \cite{LMW02} three of the authors defined
{\em augmented curvature-center coordinates},
(called {\em ACC-coordinates} for short), 
$\bw(C)$ of an
oriented circle $C$ with signed curvature $b$
and center $(x,y)$, which resolve this ambiguity.
 These represent $C$ by a row vector
\beql{ACC}
\bw(C) := (\bar{b} ,  b , b x_1, bx_2) 
\eeq
in which $\bar{b}$ is the signed curvature of the oriented circle
obtained by inversion in the unit circle.
The operation of inversion in the unit circle acts on $\RR^2$ by
$$
(x, y) \mapsto (x^{'}, y^{'}) = ( \frac{x}{x^2 + y^2}, \frac{y}{x^2 + y^2}).
$$
It maps a circle $C$ of finite (oriented) radius $r$ to the circle 
$\bar{C}$ with center $\bar{\bx} =
{\bx}/(|\bx|^2 - r^2 )$ and oriented radius 
$\bar{r} = \frac{r}{|\bx|^2 -r^2}$, 
having  oriented curvature  $\bar{b}= \frac{|\bx|^2 -r^2}{r}$.
If $C$ is a straight line with specified normal
direction, we determine $\bar{C}$ as the
image of this line, with orientation coming from the 
specified normal. In the degenerate case that $\bar{C}$ is a 
straight line we define
$\bar{b}=0$. In all cases,
\beql{303f}
(b x, b y) =  \frac {\bx}{r} = \frac{\bar{\bx}}{\bar{r}} = 
(\bar{b} x', \bar{b} y'),
\eeq
so that $C$ and $\bar{C}$ have the same curvature$\times$center data.

Augmented curvature-center coordinates provide a global coordinate
system:  no two distinct oriented circles have the same coordinates.
The only case to resolve is when $C$ is a straight line, i.e., $b=0$.
The relation \eqn{303f} shows that $(\bar{b}, b x,b y)$
are the curvature-center coordinates of $\bar{C}$, and if $\bar{b}
\neq 0$, this uniquely determines $\bar{C}$; inversion in the unit
circle then determines $C$. In fact, $\bar{b}$ in this case is twice
the distance of $C$ from the origin. In the remaining case, $b = \bar{b} =0$
and $C = \bar{C}$ is the unique line passing through the origin
whose unit normal is given by the remaining coordinates.

Given a collection $(C_1, C_2, C_3, C_4)$ of four oriented
circles (possibly lines) in $\rr^2$, the {\em augmented matrix}
$\bW$ associated with it is the $4 \times 4$ matrix whose
$j$-th row has entries given by the augmented curvature-center
coordinates $\bw(C_j)$ of the $j$-th circle.
The following result characterizes oriented
Descartes configurations.

%
%
%
%

\begin{theorem}[Augmented Euclidean Descartes Theorem]~\label{th32}
Given an ordered configuration $\sD$ of four oriented circles (or lines) 
$\{C_i:  1 \leq i \leq 4\}$
with 
curvatures $(b_1,~b_2,~b_3,~b_4)$ and centers
$\{(x_i, y_i)~: 1 \leq i \leq 4\}$,
let $\bW_\sD$ be the $4 \times 4$ matrix
\beql{308}
\bW_\sD :=
        \left[ \begin{array}{rrrr}
         \bar{b}_1  &  b_1 & b_1x_1 & b_1y_1 \\
         \bar{b}_2  &  b_2 & b_2x_2 & b_2y_2 \\
         \bar{b}_3  &  b_3 & b_3x_3 & b_3y_3 \\
         \bar{b}_4  &  b_4 & b_4x_4 & b_4y_4
        \end{array} \right] .
\eeq
If $\sD$ is an (ordered, oriented)  
Descartes configuration then $\bW_{\sD}$ satisfies 
\beql{307}
        \bW^T \bQ_{D}  \bW = 
        \left[ \begin{array}{rrrr}
           0 & -4 & 0  & 0 \\
          -4 &  0 & 0  & 0 \\
           0 &  0 & 2  & 0 \\
           0 &  0 & 0  & 2 \\ \end{array} \right].     
\eeq
Conversely, any real solution $\bW$ to \eqn{307} is the augmented
matrix $\bW_{\sD}$  of a unique ordered,
oriented Descartes configuration $\sD$.  
\end{theorem}

\paragraph{Proof.} This is proved as the
two-dimensional case of  Theorem 3.3  in \cite{LMW02}. $~~~\bsq$ \\

We call the quadratic form $Q_W$ defined by the matrix 
\beql{308b}
\bQ_{W}= \left[ \begin{array}{rrrr}
           0 & -4 & 0  & 0 \\
          -4 &  0 & 0  & 0 \\
           0 &  0 & 2  & 0 \\
           0 &  0 & 0  & 2 \\ \end{array} \right]      
\eeq
the {\em Wilker quadratic form}.
We name this quadratic form after  J. B. Wilker ~\cite{Wi81}, who
introduced in spherical geometry a coordinate system
analogous to augmented curvature-center coordinates, see
\cite[\S2 p. 388-390 and \S9]{Wi81}. However Wilker did not
explicitly formulate any result  exhibiting the quadratic form
$Q_{W}$, see  \cite[p. 349 ]{LMW02}. 

Theorem~\ref{th32}
identifies the set of all ordered, oriented
Descartes configurations $\DD$  with the set $\sM_{\DD}$ of real   
solutions $\bW = \bW_{\sD}$  to the matrix equation \eqn{307}.
This equation 
states that the  augmented matrix coordinates
of an oriented Descartes configuration give an intertwining map
between the Descartes form and the Wilker form.
The set $\sM_{\DD}$ has the structure of a six-dimensional
affine real-algebraic variety.

It is possible to refine the parameter space
to  a {\em moduli space} 
$\tilde{\sM}_{\DD}$
of (unordered, unoriented) Descartes configurations as an
orbifold $\sM_{\DD}/\sim$ obtained by quotienting by
a finite group of order $48$ (acting on the left). This group
is generated by the $4\times 4$ permutation
matrices (permuting rows) and $-\bI$, which reverses total
orientation.  This orbifold has singular points at those
ordered, oriented Descartes configurations which remain  
invariant under a non-trivial permutation matrix.
For our purposes it is more convenient to use 
the parameter space $\sM_{\DD}$ which is a smooth manifold.

Both the Descartes quadratic form and the Wilker quadratic form
are equivalent over the real numbers to the Lorentz
quadratic form 
$Q_{\sL}(x) := - x_0^2 + x_1^2 + x_2^3 + x_{3}^2,$
with associated matrix
\beql{308c}
\bQ_{\sL}= \left[
\begin{array}{crrr}
-1 & 0 & 0 & 0 \\
 0 & 1 & 0 & 0 \\
 0 & 0 & 1 & 0 \\
 0 & 0 & 0 & 1
\end{array}
\right].
\eeq
For any quadratic form $Q$, let  $Aut(Q)$  be the group of automorphisms
under the ``congruence action'', defined by 
$$
Aut(Q)=\{ \bU  \in GL(4, \RR) : \bU^T \bQ \bU = \bQ \},
$$
where $\bQ$ is the symmetric matrix that represents $Q$. 
The Lorentz quadratic form has a large group of automorphisms
$Aut(Q_{\sL})$, 
which is exactly the real Lorentz group $O(3, 1)$.

The Descartes form and Wilker form are  not only equivalent
to the Lorentz form  over the real numbers,
but also  over the rational numbers. For the Descartes
form one  has
\beql{309aa}
\bQ_{D} = \bJ_0^T \bQ_{\sL} \bJ_0,
\eeq
in which
\begin{equation}~\label{309}
\bJ_0 = \frac{1}{2} \left[
\begin{array}{crrr}
1 & 1 & 1 & 1 \\
1 & 1 & -1 & -1 \\
1 & -1 & 1 & -1 \\
1 & -1 & -1 & 1
\end{array}
\right]
\end{equation}
and $\bJ_0 = \bJ_0^T = \bJ_0^{-1}$.
It follows that
\beql{N309c}
Aut(Q_{D})=\bJ_{0}^{-1} O(3,1) \bJ_0.
\eeq

\begin{figure}[htbp]
\centerline{\epsfxsize=3.0in \epsfbox{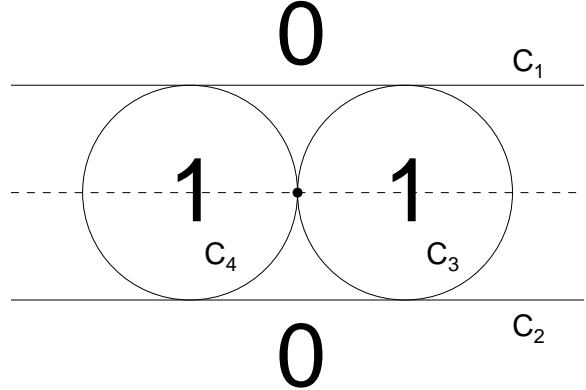}}
\caption{A strongly integral Descartes configuration}~\label{fig2.1a}
\end{figure}

The rational equivalence
of the Wilker quadratic form $Q_{W}$ 
to the Lorentz form $Q_{\sL}$ follows from Theorem~\ref{th32},
as soon as we exhibit a  Descartes configuration $\sD$
whose augmented curvature-center coordinates $\bW_{\sD}$ are
a rational matrix. One is given by
\begin{equation}~\label{N330a}
\bW_0= \bW_{\sD}= \left[
\begin{array}{rrrr}
 2 & 0 & 0 & 1 \\
 2 & 0 & 0 & -1 \\
 0 & 1 & 1 & 0 \\
 0 & 1 & -1 & 0
\end{array}
\right].
\end{equation}
It corresponds to the positively oriented Descartes configuration pictured in
Figure \ref{fig2.1a}, in which the dotted line is the $x$-axis
and the two circles touch at the origin.
%
%
%
%

We then have 
\begin{equation}~\label{N333a}
\bQ_{W} = \bA^{T} \bQ_{\sL} \bA,  
\end{equation}
with 
\begin{equation}~\label{N333b}
\bA =  \bJ_0 \bW_0=\left[
\begin{array}{crrr}
 2 & 1 &  0 & 0 \\
 2 & -1 & 0 & 0 \\
 0 & 0 & 1 &  1 \\
 0 & 0 & -1 & 1
\end{array}
\right]. 
\end{equation}
It follows that
\beql{N333c}
Aut(Q_W)=\bA^{-1} O(3,1) \bA.
\eeq

Finally we note that the determinants of these quadratic
forms are given by  
\beql{Ndet1}
\det( \bQ_{D}) = \det(\bQ_{\sL}) = -1~~\mbox{and}~~ \det(\bQ_{W}) = -64.
\eeq
Taking determinants in  \eqn{307} yields 
that augmented curvature-center coordinate
matrices $\bW_{\sD}$ have
\beql{Ndet2}
\det( \bW_{\sD}) = \pm 8. 
\eeq
Examples show that  both  values $\pm 8$ occur.

%
%
\vspace{.3cm} 
\noindent {\bf Proof of the ``only if'' part of Theorem \ref{th31}.} 
We note that if $\bM$ is the curvature-center coordinate
matrix of a Descartes configuration, this configuration
is unique. Indeed at least two circles in the configuration
have nonzero curvature, so are determined by their
curvature-center coordinates. This gives enough information
to find the location of any lines in the configuration,
and their signed normal vectors  are determined uniquely
by their curvature-center coordinates.
Thus it  suffices to show that there exists a $4 \times 4$ matrix $\bW$
whose last three columns agree with $\bM$, such that
\begin{equation}~\label{A1A}
\bW^T \bQ_{D} \bW = \bQ_{W} = \left[ \begin{array}{rrrr}
           0 & -4 & 0  & 0 \\
          -4 &  0 & 0  & 0 \\
           0 &  0 & 2  & 0 \\
           0 &  0 & 0  & 2 \\ \end{array} \right].
\end{equation}
If so,  then Theorem \ref{th32}  implies
that $\bW$ is the augmented curvature-center coordinate
matrix of a unique ordered, oriented Descartes configuration $\sD$.
It follows that $\bM$ is the curvature-center coordinate matrix
of $\sD$.

This discussion implies that the extension $\bW$ is unique,
and we proceed to find it by transforming the problem
from the Descartes form $\bQ_{D}$ to the Lorentz form
$\bQ_{\sL} = \bJ_0^T \bQ_{D} \bJ_0$.  This proof
seems ``pulled out of thin air'' in that it  guesses to do
the ``translation'' by the matrix $\bT$ below, in terms of
which the extended matrix $\tilde{\bW}$ has a simple form.
We set
$\tilde{\bM} = \bJ_0 \bM$
with $\bJ_0$ given by \eqn{309}, so that  $\bJ_0 = \bJ_0^T = \bJ_0^{-1}$.
Then
\beql{A2}
\tilde{\bM}^T \bQ_{\sL}\tilde{\bM} = \bM^T \bQ_{D} \bM =  \bQ_{0} :=
 \left[
\begin{array}{ccc}
           0 & 0 & 0 \\
           0 & 2 & 0 \\
           0 & 0 & 2
        \end{array} \right].
\eeq
We now note that
$$
\tilde{\bM}_{11} = \frac{1}{2}( m_{11} + m_{12} + m_{13} + m_{14}) \ne 0.
$$
To verify this, note that the top left entry of  \eqn{209} gives
$$
Q_{D}(m_{11},  m_{12}, m_{13}, m_{14}) =
(m_{11}^2 + m_{12}^2 + m_{13}^2+ m_{14}^2 ) -
\frac{1}{2} ( m_{11} + m_{12} + m_{13} + m_{14})^2 = 0,
$$
and if $m_{11} + m_{12} + m_{13} + m_{14}= 0$ then we would obtain
$m_{11}^2 + m_{12}^2 + m_{13}^2+ m_{14}^2=0$, which forces
all  $m_{1j}=0$, contradicting the hypothesis.
Now there exists a unique matrix
$$
\bT= \left[
\begin{array}{ccc}
 a & b & c \\
 0 & 1 & 0 \\
 0 & 0 & 1
\end{array} \right]
$$
such that the first row of $\tilde{\bM} \bT$ is $[ 1, 0, 0]$,
which takes $a = (\tilde{\bM}_{11})^{-1}$,  so that $\bT$ is invertible.
We write
$$
\tilde{\bM} \bT = \left[
\begin{array}{ccc}
 1 & 0 & 0 \\
 v_{21} & v_{22} \sqrt{2} & v_{23} \sqrt{2} \\
 v_{31} & v_{32} \sqrt{2} & v_{33} \sqrt{2} \\
 v_{41} & v_{42} \sqrt{2} & v_{43} \sqrt{2}
\end{array}
\right] ~,
$$
and then have
\begin{equation}\label{A3}
\bT^T \tilde{\bM}^T \bQ_{\sL} \tilde{\bM} \bT =
\bT^T  \bQ_0 \bT = \bQ_0 = \left[
       \begin{array}{ccc}
           0 & 0 & 0 \\
           0 & 2 & 0 \\
           0 & 0 & 2
        \end{array}
\right].
\end{equation}
This  matrix  equation is equivalent to the assertion that
\begin{equation*}
\bV := \left[
\begin{array}{ccc}
v_{21} & v_{22} & v_{23} \\
v_{31} & v_{32} & v_{33} \\
v_{41} & v_{42} & v_{43}
\end{array} \right]
\end{equation*}
is an orthogonal matrix, i.e. $\bV^T \bV = \bV \bV^T =I.$
We now define
$$
\tilde{\bW} :=
\left[
\begin{array}{cccc}
2 &  1 & 0 & 0 \\
-2 v_{21} &  v_{21} & v_{22} \sqrt{2} & v_{23} \sqrt{2} \\
-2 v_{31} &  v_{31} & v_{32} \sqrt{2} & v_{33} \sqrt{2} \\
-2 v_{41} &  v_{41} & v_{42} \sqrt{2} & v_{43} \sqrt{2}
\end{array}
\right],
$$
whose last three columns match $\tilde{\bM} \bT$,
and verify by direct calculation using \eqn{A3}  that
$$
\tilde{\bW}^T \bQ_{\sL} \tilde{\bW} = \bQ_{W}.
$$
We define $\bW := \bJ_0^{-1} \tilde{\bW} = \bJ_0 \tilde{\bW}$,
and find that
$$
\bW^T \bQ_{D} \bW = \tilde{\bW}^T \bJ_0  \bQ_{\sL} \bJ_0 \tilde{\bW} = \bQ_{W}.
$$
This is the desired lift,
since  the last three columns of $\bW$ are exactly $\bM$.$~~~\bsq$

%
%

\subsection{M\"{o}bius and Lorentz group actions}

The augmented Euclidean Descartes theorem immediately 
yields two group actions on the space $\sM_{\DD}$
of ordered, oriented Descartes configurations.
The group $Aut( Q_{D})$ acts on the left
and the group $Aut(Q_{W})$ acts on the right, as
$$
\bW_{\sD}  \mapsto \bU  \bW_{\sD}  \bV^{-1}, ~~\mbox{with}~~ 
\bU \in Aut( Q_{D}),~ 
\bV \in Aut(Q_{W}).
$$
The two group actions clearly commute with each other.
Both groups are  conjugate to the real Lorentz group
$O(3, 1)$, and therefore they each have four connected 
components. These components are specified for any  $\bY \in O(3,1)$
by the
sign of $\det(\bY)$,  and by the sign of ``total orientation'', 
which is the sign of  
$\bY_{11} = \be^T_1 \bY \be_1$, in which $\be^T_1=(1,0,0,0)$. 
For  $\bU \in Aut(Q_{D})$, it is the sign of $\bo^T \bU \bo$, in which 
$\bo^T=(1,1,1,1)$; and for $\bV\in Aut(Q_W)$, it is the sign of 
$\be^T_1 \bA \bV \bA^{-1} \be_1 $ in which  $\bA$ is any matrix satisfying 
$O(3,1)=\bA  Aut(Q_W)\bA^{-1}$, such as \eqn{N333b}.  

The parameter space $\sM_{\DD}$ also has four connected
components, specified by similar invariants, which are the sign
of $\det(\bW)$ and the {\em total orientation}, which is the sign of the
sum\footnote{The sum of the signed curvatures of the circles 
of an oriented  Descartes configuration cannot
be zero.}  
of the (signed) curvatures of the four circles in the
Descartes configuration.
We let
$$
\sM_{\DD} = \sM_{+}^{\uparrow} \cup \sM_{-}^{\uparrow} \cup 
\sM_{+}^{\downarrow} \cup \sM_{-}^{\downarrow}, 
$$
in which the subscript describes the sign of the determinant and
the superscript the orientation, with $\uparrow$ being 
positive orientation. We also let
$$
\sM_{\DD}^{\uparrow} := \sM_{+}^{\uparrow} \cup \sM_{-}^{\uparrow} 
$$
denote the set of positively oriented Descartes configurations.
The relevance of this decomposition is that the Apollonian group
defined in \S4 leaves the set of positively oriented Descartes
configurations 
$\sM_{\DD}^{\uparrow}$ invariant.

The action on the right by the elements of $Aut(Q_{W})$
maps circles to circles, since the 
circles in a Descartes configuration correspond 
to the rows in the matrix $\bW_{\sD}$ of an ordered, oriented
Descartes configuration.  
This right action 
can essentially be identified
with the general M\"{o}bius group M\"{o}b(2) generated
by the  linear fractional
transformations acting on the one-point compactification
$\hat{\RR}^2$ of $\RR^2$, which is $PSL(2, \CC) = SL(2, \CC)/\{ \pm \bI\}$,
together with  complex conjugation $\bz \mapsto \bar{\bz}$.
and has two connected components.
More precisely,
$Aut(Q_W)$ is isomorphic to 
the direct sum of the M\"{o}bius group with a group of
order two, which has four connected components,
as explained in Appendix A.

The action on the left, by $Aut(Q_{D})$,
produces a new oriented Descartes configuration
whose coordinates  mix together the coordinates of the 
different circles in the original   Descartes configuration.
This group action  does
not make sense as an action on individual circles
in the configuration.
This group action is intrinsically associated to
the $6$-dimensional space of oriented
Descartes configurations.

\begin{theorem}~\label{th33}
(1) The groups $Aut(Q_{D})$ and $Aut(Q_{W})$ are
conjugate to $Aut(Q_{\sL}) \equiv  O(3,1)$.

(2) The group $Aut(Q_{D})$ acts transitively on the 
left on the space $\sM_{\DD}$
of all ordered, oriented  Descartes configurations. 
Given two such 
 Descartes configurations $\sD$ and $\sD^{'}$
there exists a unique $\bU \in Aut(Q_{D})$ such that
$\bU \bW_{\sD} = \bW_{\sD'}$.

(3) The group $Aut(Q_{W})$ acts transitively on the 
right on the space 
of all ordered, oriented  Descartes configurations $\sM_{\DD}$. 
Given two such Descartes configurations $\sD$ and $\sD^{'}$
there exists a unique $\bV \in Aut(Q_{W})$ such that
$\bW_{\sD} \bV^{-1} = \bW_{\sD'}$.

(4) The action of $Aut(Q_{D})$ on the space 
$\sM_{\DD}$ commutes with  the action of $Aut(Q_{W})$. 
\end{theorem}

\paragraph{Remark.} The left action by $Aut(Q_{D})$
and the right action by $Aut(Q_W)$ on $\sM_{\DD}$
given by Theorem \ref{th33} can be identified with a left and right
action of the Lorentz group $O(3,1) = Aut(Q_{\sL})$ on $\sM_{\DD},$
using \eqn{N309c} and \eqn{N333c}, respectively.
Theorem~\ref{th33} 
shows that both these actions are  transitive, and that  the
stabilizer of a point is the identity element.
This is equivalent to saying that 
the space $\sM_{\DD}$ is a {\em principal homogeneous space} 
(or {\em torsor}) for 
$O(3,1)$ for either action.

\paragraph{Proof.}
Part (1) follows from the conjugacy between the
Descartes and Wilker forms and the Lorentz form $Q_{\sL}$
given in \eqn{309aa} and \eqn{N333a}, respectively.
These give
$$
Aut(Q_{D}) = \bJ_0^{-1} Aut(Q_\sL) \bJ_0 \equiv \bJ_0^{-1} O(3,1) \bJ_0,
$$
and
$$
Aut(Q_{W}) = \bA^{-1} Aut(Q_{\sL}) \bA \equiv \bA^{-1} O(3,1) \bA.
$$

Parts (2) and (3) follow immediately from (1). Given a fixed
$\bW= \bW_{\sD} \in \sM_{\dd}$, we assert that
\begin{equation}~\label{N335}
\sM_{\DD} = Aut(Q_{D}) \bW. 
\end{equation}
Taking  $\bW_0$ in \eqn{N330a} we have
$$
(\bW \bW_0^{-1})^T \bQ_{D} (\bW \bW_0^{-1})=
(\bW \bA^{-1}\bJ_0)^{T} \bQ_{D} (\bW \bA^{-1}\bJ_0) = \bQ_{D},
$$
so that $\bW \bW_0^{-1} \in Aut(Q_{D})$ and
 $Aut(Q_{D}) \bW \bW_0^{-1}$ forms  a single orbit of $Aut(Q_{D})$,
The map $\bW_{\sD} \mapsto \bW_{\sD} \bW_0^{-1}$ from 
$\sM_{\dd}$ into $Aut(Q_{D})$ is one-to-one since $\bW_0$
is invertible, and it is onto since the domain includes
$Aut(Q_{D}) \bW$. Thus 
\eqn{N335} follows, and this 
gives (2).  We similarly obtain 
$$
\sM_{\DD} = \bW Aut(Q_{W})^{-1} = \bW Aut(Q_{W}),
$$
which gives (3).

Finally, part (4) follows directly from parts (2) and (3),
since in general a left and a right action of two matrix groups
on a space of matrices commute.
~~~$\bsq$ \\

In this paper we mainly study  structures 
associated to the left action of
$Aut(Q_{D})$ on $\sM_{\DD}$.  We term the left action the 
{\em Lorentz action}, although this  is a misnomer, since the right
action can also be identified with a Lorentz group action. 
However we wish to assign different names for the two actions, 
and the right
action is conveniently called
the {\em M\"{o}bius action},
because of its relation to a 
M\"{o}bius group action on circles on the Riemann sphere,
detailed in Appendix A.

%
%

\subsection{Distinguished Elements of $Aut(Q_{D})$.}

We now describe 
some  specific elements of $Aut(Q_{D})$ that have a 
nice geometrically
interpretable action on every Descartes configuration,
visualizable in terms of inversions, and whose 
associated matrices have integer (or half-integer) entries.

\begin{figure}[htbp]
\centerline{\epsfxsize=3.0in \epsfbox{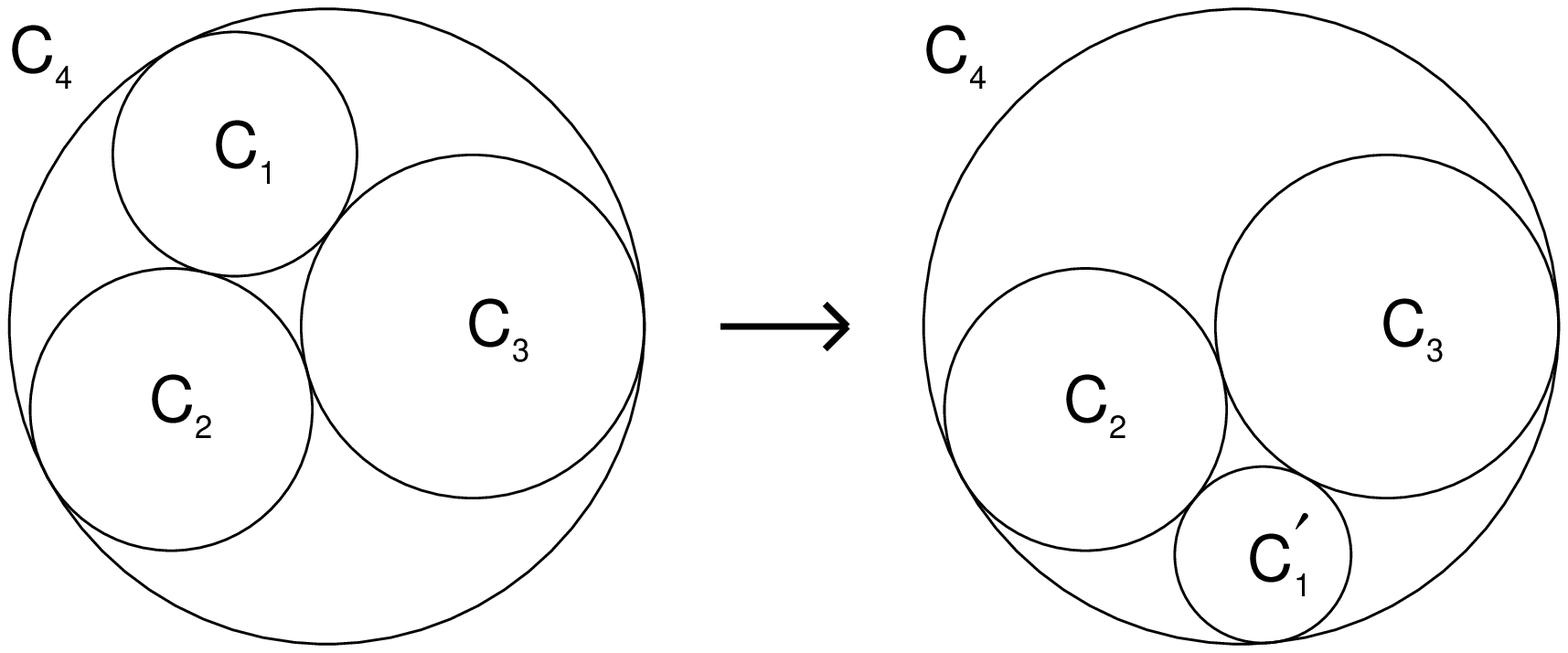}}
\caption{The reflection operator $\fs_1$ }~\label{fig4}
\end{figure}

The first set of four operations correspond to inversion in 
the circle determined by the three intersection points  of circles
in a Descartes configuration $\sD$ that avoid  one 
particular circle. There are four possibilities for these.
This inversion fixes the three circles involved in
the intersections and moves the fourth circle, to
the unique other circle that is tangent to the first
three circles. 
For this reason we call it the  {\em reflection operator}. 
Let 
$\fs_1=\fs_1[\sD]$ denote the M\"{o}bius
transformation  of this kind that moves the 
circle $C_1$. In particular, $\fs_1$ maps the Descartes 
configuration $\sD=(C_1, C_2, C_3, C_4)$ to $\fs_1(\sD)=(C_1', C_2, 
C_3, C_4)$.  
This reflection operator depends on the specific Descartes configuration.
However for  all Descartes configurations $\sD$ there holds 
\beql{N331}
  \bW_{\fs_1(\sD)}=  \bS_1 \bW_{\sD}  ,
\eeq
where  $\bS_1 \in Aut(Q_{D})$ is given by 
$$
\bS_1  =  \left[
\begin{array}{rccc}
-1 & 2 & 2 & 2 \\
0 & 1 & 0 & 0 \\ 
0 & 0 & 1 & 0 \\
0 & 0 & 0 & 1
\end{array} \right]
$$
belongs to $Aut(Q_{D})$ and is independent of the location of $\sD$.
The geometric action of $\fs_1$  is pictured in Figure~\ref{fig4}.
The other three operations $\fs_2, \fs_3, \fs_4$ give similar matrices
$\bS_2, \bS_3, \bS_4$, 
obtained by permuting the first and j-th rows and columns of $\bS_1$
for $2 \le j \le 4$.

%
%
%

To prove \eqn{N331} we check  that 
$\bS_1 \in Aut(Q_{D})$. 
Then for any Descartes
configuration $\sD$  Theorem~\ref{th32} gives 
that  $\bS_1 \bW_{\sD} =\bW_{\sD'}$ for some
Descartes configuration $\sD'$.  The configuration
$\sD'$ necessarily has three oriented circles fixed, with the fourth
moved, so the fourth one must be the unique
other circle tangent to the given three, and its orientation
is uniquely determined by the other three orientations. Now
$\bW_{\fs_1(\sD)}$ corresponds to  that Descartes
configuration  consisting  of the same four ordered circles,
which has the same (positive or negative) orientation as $\sD'$. 

\begin{figure}[htbp]
\centerline{\epsfxsize=3.0in \epsfbox{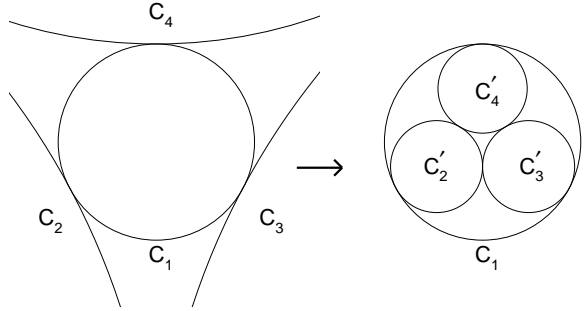}}
\caption{The inversion operation $\fs_1^\perp$ }~\label{fig5}
\end{figure}

A second set of four operations corresponds to inversion in
one of the four circles in a Descartes configuration.
Let that circle be $C_1$. Now $C_1$ remains fixed,
while the other three circles change. Denoting
this inversion by $\fs_1^{\perp}$,  it is easy to show that
\beql{N323a}
\bW_{\fs_1^{\perp}(\sD)} = \bS_1^{\perp} \bW_{\sD},
\eeq
where $\bS_1^{\perp} \in Aut(Q_{D})$ is given by
$$
\bS_1^{\perp} =  \left[
\begin{array}{rccc}
-1 & 0 & 0 & 0 \\
2 & 1 & 0 & 0 \\
2 & 0 & 1 & 0 \\
2 & 0 & 0 & 1
\end{array} \right].
$$
The corresponding matrices $\bS_2^{\perp},\bS_3^{\perp},\bS_4^{\perp}$
are obtained by permuting the first and j-th rows and columns of
$\bS_1^{\perp}$ for $2 \le j \le 4$.
The operation $\fs_1^{\perp}$ is pictured in Figure~\ref{fig5}.

%
%

\begin{figure}[htbp]
\centerline{\epsfxsize=2.0in \epsfbox{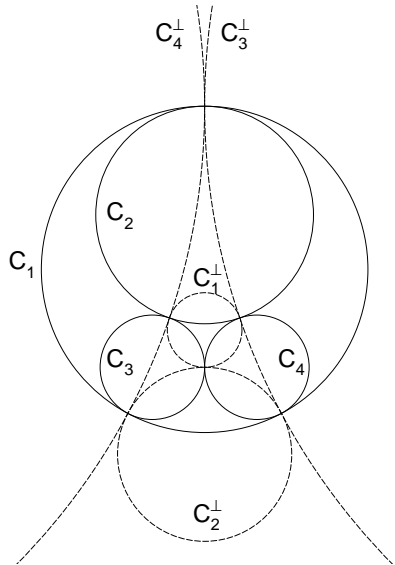}}
\caption{The dual operation}~\label{dualfig}
\end{figure}

Finally we describe an operator, 
which we call the {\em duality operator} $\bD$,  
and which is an involution, as follows.
Given a Descartes configuration $\sD$, there exists a
Descartes configuration $\sD'$ obtained from it which
consists of the four circles each of which passes through
the three intersection points avoiding one circle.
The resulting configuration has the same six points of
tangency as the original configuration, and the circles
in the configuration $\sD'$ are perpendicular to the
circles of $\sD$ at these tangency points. This is
pictured in Figure~\ref{dualfig}.

%
%

Straightforward  computation shows that
$$
\bW_{\sD'}  = \bD \bW_{\sD},
$$
with
$$
\bD = \frac{1}{2} \left[
\begin{array}{crrr}
-1 & 1 & 1 &  1 \\
1 & -1 &  1 &  1 \\
1 &  1 & -1 &  1 \\
1 & 1 &  1 & -1
\end{array}
\right].
$$
Note that  $\bD = -\bQ_{D} \in Aut (Q_{D})$. There is
a M\"{o}bius transformation $\fd = \fd_{\sD}$ which
sends $\sD \mapsto \fd(\sD) =  \sD'$, and  it
depends on $\sD$.

%
%
%

\section{Apollonian Packings and the Apollonian Group}
\setcounter{equation}{0}

\begin{figure}[htbp]
\centerline{\epsfxsize=2.0in \epsfbox{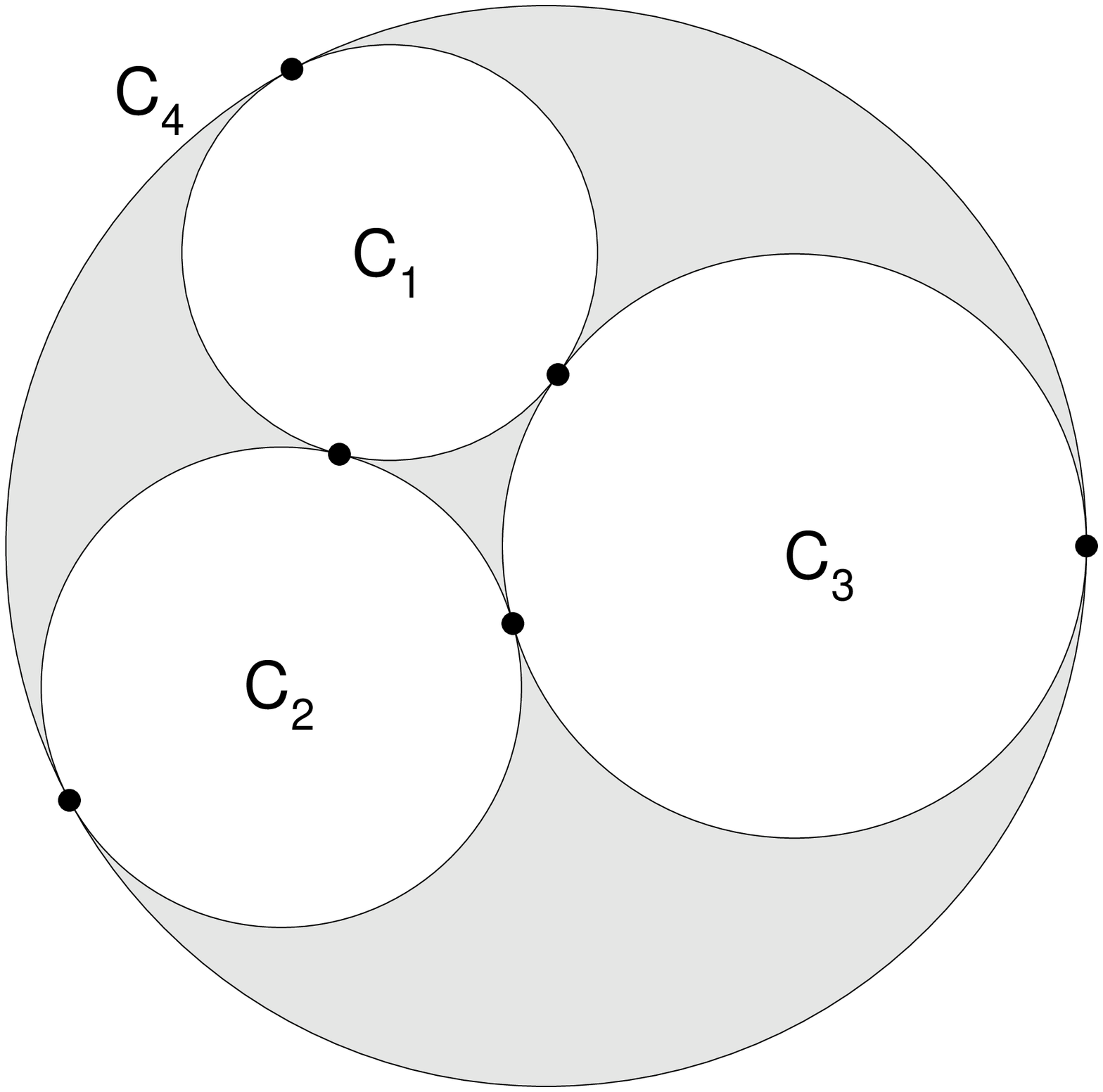}}
\caption{A zero-th stage packing $\sP_{\sD}^{(0)}$ 
and uncovered regions (lunes)}\label{fig3.1}
\end{figure}

Apollonian circle packings are infinite packings of circles
recursively constructed from a given 
positively oriented Descartes configuration $\sD$.
For simplicity consider a
positively oriented  Descartes configuration
$\sD = \{C_1, C_2, C_3, C_4\}$ in which circle $C_4$ encloses the other three,
so that the interior of $C_4$ includes the point $\bz_\infty$ at infinity.
The uncovered area consists of four lunes, indicated by the shading 
in Figure~\ref{fig3.1}.

%
%

\begin{figure}[htbp]
\centerline{\epsfxsize=2.0in \epsfbox{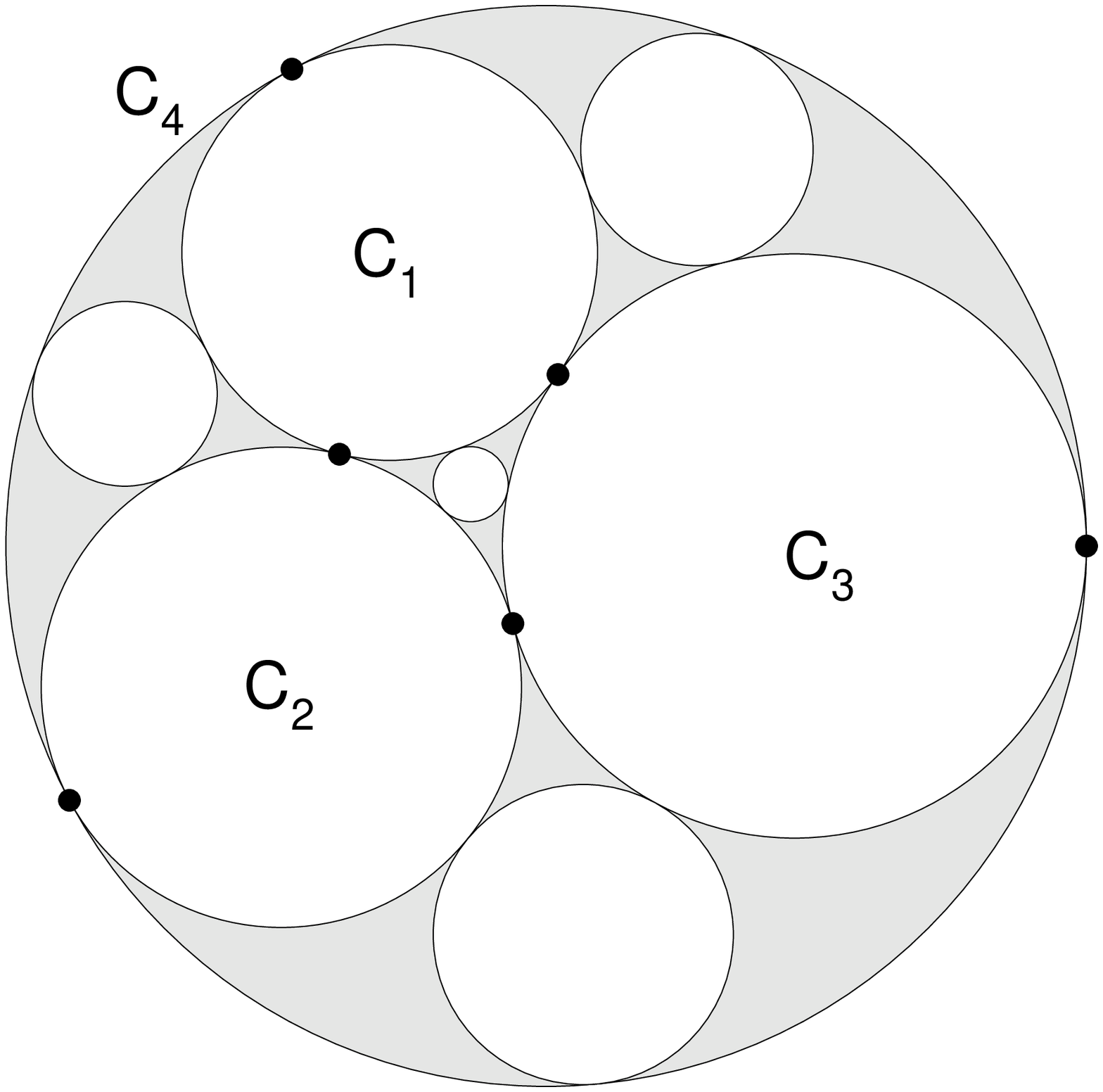}}
\caption{A first  stage packing $\sP_{\sD}^{(1)}$ }\label{fig3.2}
\end{figure}

The zeroth stage packing $\sP_{\sD}^{(0)} = \sD$.
At the first stage we inscribe a circle in each lune,
to obtain a circle packing $\sP_{\sD}^{(1)}$ containing 8 circles.
Each of these circles lies in a unique (unordered) Descartes configuration
in the first stage packing, consisting of it and the three circles it touches.
The uncovered area that remains consists of 12 lunes. See Figure \ref{fig3.2}.

%
%
 
At the second stage we inscribe a circle in each of these lunes, which 
produces the second stage packing $\sP_{\sD}^{(2)}$.
Continuing in this way, at the $n$-th stage we add $4 \cdot 3^{n-1}$ circles, 
and each of these lies in a unique
Descartes configuration in the $n$-th stage partial packing.
The {\em Apollonian packing} $\sP_{\sD}$ associated to $\sD$ is
the limit packing 
\beql{T31}
\sP_{\sD} := \bigcup_{n=1}^\infty \PD^{(n)} ~.
\eeq

One can  regard an Apollonian packing $\sP$ 
as a geometric object consisting of an infinite collection of
circles. These circles are described  as the four
orbits 
of a group of M\"{o}bius transformations $G_{\sA}(\sD)$
acting on the circles
$(C_1, C_2, C_3, C_4)$ in
original Descartes configuration $\sD$. The group
\beql{4311}
G_{\sA}(\sD) := \langle \fs_1, \fs_2, \fs_3, \fs_4\rangle, 
\eeq
in
which  $\fs_j$ is the inversion with respect to the circle that
passes through the three intersection points of $\sD$ that do not
include the circle $C_j$, as indicated in \S3.3. One can check
that the circles added at the $n$-th stage of the construction above
correspond to words $\fs_{i_1} \fs_{i_2} \cdots \fs_{i_n}$
of length $n$ in the generators of this group,
in which any two adjacent generators are distinct, i.e. $i_j \ne i_{j+1}$.

The key geometric property of Apollonian circle packings
is that they can be viewed as packings of disks having
the circles as boundary. 
To each circle in the Riemann sphere correspond
two  disks having the
circle as boundary. Recall that the ``interior'' of an
oriented circle corresponds to
making  a choice of one of these two disks, 
indicated by a choice of sign of the curvature, when it is
nonzero (and by a choice of normal vector if 
the curvature is zero). 
The disk packing property 
of an Apollonian packing  corresponds to treating its 
Descartes configurations as positively oriented.  

%
%

\begin{theorem}~\label{th40}
For any  Apollonian packing generated by a
positively oriented Descartes configuration,
the interiors of all circles in the packing
are disjoint.
\end{theorem}

\paragraph{Proof.} 
We generate the packing from a single
Descartes configuration $\sD$ chosen to have
positive total orientation, so that the interiors
of its four circles defined by the total orientation
are disjoint. 
Positive total  orientation  is preserved by the action of M\"{o}bius
transformations; see Theorem~\ref{th90a} in Appendix A.
Let 
$$
\sP_{\sD}[m]:=\bigcup_{n=1}^m \sP_{\sD}^{(n)},
$$
where $\sD$ is a positively oriented Descartes configuration. 
We prove by induction on $m$ the following statement: 
 $\sP_{\sD}[m]$ consists of $2(3^m+1)$ circles which form the boundary 
of $2(3^m+1)$ disks and $4\cdot 3^m$ lune areas. The  interiors of the 
disks, as well as the lunes, are mutually disjoint. 

The base case $m=0$ is easy. Since any two positively oriented Descartes 
configurations are equivalent under the M\"{o}bius transformations, which
map circles to circles and preserve the (total) orientation, 
(c.f. Theorem 7.2, Appendix A),
we can simply check the Descartes configuration given in Figure 7, which 
has $4$ circles that form the boundary of $4$ disks and $4$ lunes 
with disjoint 
interiors. 

Assume the statement holds for $m \geq 0$. We carry out 
 the inductive step for  $m+1$. By the inductive construction above, 
\[
\sP_{\sD}[m+1] = \sP_{\sD}[m] \cup \sP_{\sD}^{(m+1)}.
\]
The circles of $\sP_{\sD}[m]$ remain in $\sP_{\sD}[m+1]$ with their
interiors untouched. In each lune of $\sP_{\sD}[m]$ a new circle is
inscribed, breaking the original lune into 4 parts--one disk and three
smaller lunes, with disjoint interiors. Combining the inductive hypothesis,
we have that  
$\sP_{\sD}[m+1]$ consists of $2(3^m+1)+4\cdot 3^m =2(3^{m+1}+1)$ many 
circles with empty interiors, and the uncovered area are $3\cdot 4\cdot 3^m
=4\cdot 3^{m+1}$ lunes  with disjoint interiors.  
Each new Descartes configuration shares the interior of three circles
with the Descartes configuration  generating it; this implies it has
the same orientation as the preceding one; hence it has positive
total orientation.
This proves the 
statement for $m+1$.  $~~~\bsq$ \\

We first consider an Apollonian packing as a geometric
object. A {\em geometric Apollonian packing} $\sP$
is the point set
consisting of a countable collection  of circles on the
Riemann sphere $\hat{\cc} = \cc \cup \{\infty\}$
obtained by 
the construction above. It has Hausdorff dimension $1$,
since it is a countable union of sets of Hausdorff dimension $1$
(circles).
Let  $G(\sP)$ denote the group of M\"{o}bius
transformations that leave $\sP$ invariant. This group can be proved
to be a discrete subgroup of the group of all M\"{o}bius transformations,
which acts transitively on the circles in the packing.
The group $G(\sP)$ contains 
$G_{\sA}(\sD)$ above as a subgroup of index $24$,
with cosets given by 
$24$ M\"{o}bius transformations whose effect is to
fix a generating Descartes configuration $\sD$ and permute the
four circles in it.
(These are nontrivial facts, and we do not prove them here.)

The geometric Apollonian packing is not a closed set on the
Riemann sphere. We define  $\Lambda(\sP)$ to be its closure on
the Riemann sphere, and call it the {\em  residual set}
of the Apollonian packing. These sets are prototypical examples
of fractal sets, and  
have been much studied; in 1967 Hirst \cite{Hi67}
showed these sets have Hausdorff dimension strictly between
$1$ and $2$. For later use, we
summarize properties of  $\Lambda(\sP)$   in the following result. 
%
%
%
%

\begin{theorem}~\label{th41b}
The residual set $\Lambda(\sP)$ of a
geometric  Apollonian packing $\sP$ has the
following properties.

(1) $\Lambda(\sP)$  is the complement in
the Riemann sphere of the interiors of
all circles in the packing. Here  ``interior'' is defined
by a positive orientation of a Descartes configuration in
the packing.

(2) $\Lambda(\sP)$ has a Hausdorff dimension
which is independent of the packing $\sP$, that satisfies
the bounds 
$$
   1.300197 < \mbox{dim}_{H}(\Lambda(\sP)) < 1.314534.
$$
In particular, $\Lambda(\sP)$   has Lebesgue measure zero.

(3) $\Lambda(\sP)$ is the closure of the countable set of
all tangency points of circles in the packing.

(4) $\Lambda(\sP)$ is invariant under the action of the
discrete group $G(\sP)$ of M\"{o}bius transformations.
\end{theorem}

\paragraph{Remark.} Property (1) is the more usual
definition of residual set of an Apollonian packing.

\paragraph{Proof.}
Let $I(\sP)$ be  the open set which is the union of  the
interiors of all circles in the packing. 
It is clear that 
$$
\Lambda(\sP) \subset \Lambda^{*}(\sP) := \hat{\cc} \backslash I(\sP).
$$
The set $\Lambda^{*}(\sP)$ is the usual definition of  the
residual set of an Apollonian packing; we will show
$\Lambda(\sP) =  \Lambda^{*}(\sP).$
  
The result of Hurst \cite{Hi67} that the Hausdorff
dimension of $\Lambda^{*}(\sP)$
is strictly less than $2$  implies
that $\Lambda^{*}(\sP)$
 has Lebesgue measure zero. 
The sharper bounds on the Hausdorff dimension 
stated  here were obtained in 1973 by 
Boyd ~\cite{Bo73b}. A later result of Boyd \cite{Bo82}
suggests that 
$\mbox{dim}_{H}(\Lambda^{*}(\sP)) = 1.3056 \pm .0001.$
Thomas and Dhar\cite{TD94} give a non-rigorous
approximation scheme  suggesting that this Hausdorff dimension
is $1.30568673$ with an error of $1$ in the last digit.
See Falconer \cite[pp. 125--131]{Fa85} for proofs
of $1 < \mbox{dim}_{H}(\Lambda^{*}(\sP)) < 1.432.$

Now let $\bz \in \Lambda^{*}(\sP)$. The
interiors of the
circles in the Apollonian packing cover all but a
measure zero area in a neighborhood of  $\bz$.
This requires infinitely many circles,
so their radii must go to zero, since there are only finitely
many circles of radius larger than any positive bound.
Since each of these circles has a tangency point on it,
the point $\bz$ is a limit point of such tangency points.
Since the tangency points are contained in $\sP$,
so we conclude that 
$\Lambda^{*}(\sP) \subset \Lambda(\sP).$
This gives $\Lambda(\sP) = \Lambda^{*}(\sP)$, at
which point properties (1), (2) and (3) follow. 

Property (4) follows by
observing that the group invariance of $\sP$ carries over
to its closure, by applying it to any Cauchy
sequence of points in $\sP$. $~~~\bsq$. \\

The main viewpoint of this paper is to treat
an Apollonian circle packing $\sP$
as  described ``algebraically''  by
the set of all Descartes configurations it contains,
which we denote $\DD(\sP)$. In  the construction process above, 
except for the four circles in $\sD$, each other  circle
$C$ in the packing corresponds to the unique (unordered) 
Descartes configuration
containing it which occurs at the stage of the construction where $C$ first 
appears. Each (unordered, unoriented) Descartes configuration  appears exactly
once in this correspondence,  except for the base configuration $\sD$
which corresponds to the initial  four circles.
These Descartes configurations can be described in terms of
orbits of a discrete subgroup $\sA$ of  $Aut(Q_{\sD})$
called here the Apollonian group.

The Apollonian group $\sA$ is a subgroup of $Aut(Q_{D})$
defined by the action of
the reflection operations on Descartes configurations
described in \S3.3. These correspond to inversions in the circles
passing through the three tangency
points in a Descartes configuration that do not include
one fixed circle in the configuration.

\begin{defi}\label{de42}
{\rm 
The {\em Apollonian group} $\sA$ is the 
subgroup
of $Aut(Q_{D})$ defined  by 
\beql{Y318a}
\sA : = \langle \bS_1, \bS_2, \bS_3, \bS_4 \rangle ~,
\eeq
where
\beql{Y319a}
\begin{array}{ccclccc}
\bS_1 & = & \left[
\begin{array}{rccc}
-1 & 2 & 2 & 2 \\
0 & 1 & 0 & 0 \\ 
0 & 0 & 1 & 0 \\
0 & 0 & 0 & 1
\end{array}
\right]~, &~~~~~~ & \bS_2 & = & \left[
\begin{array}{crcc}
1 & 0 & 0 & 0 \\
2 & -1 & 2 & 2 \\
0 & 0 & 0 & 0 \\
0 & 0 & 0 & 0
\end{array}
\right] ~, \\ [+.2in]
\bS_3 & = & \left[
\begin{array}{ccrc}
1 & 0 & 0 & 0 \\
0 & 1 & 0 & 0 \\
2 & 2 & -1 & 2 \\
0 & 0 & 0 & 1
\end{array}
\right] ~, &&
\bS_4 & = & \left[
\begin{array}{cccr}
1 & 0 & 0 & 0 \\
0 & 1 & 0 & 0 \\
0 & 0 & 1 & 0 \\
2 & 2 & 2 & -1
\end{array}
\right] ~.
\end{array}
\eeq
}
\end{defi}

We now characterize
the set of ordered, oriented Descartes configurations $\DD(\sP)$ 
in an Apollonian
circle packing $\sP$ in terms of this group. We note that
a single unordered, unoriented
Descartes configuration,
corresponds to  exactly 48 ordered,
oriented Descartes
configuration, since there are
 24 choices of ordering of the four circles, and
two choices of orientation.

%
%

\begin{theorem}~\label{th41}
 The ordered, oriented Descartes configurations in an Apollonian circle 
packing comprise a  union of 48  orbits of the Apollonian group.
Each of these orbits contains exactly one ordered,
oriented representative of each 
(unordered, unoriented) Descartes configuration in the packing.
\end{theorem}

\paragraph{Proof.}
Let $\sD=(C_1, C_2, C_3, C_4)$ be the zeroth stage packing of $\sP$, which is
ordered and positively oriented. It is enough to show that
$\sA(\sD)$ contains exactly one ordered, positively
oriented representative of each
unordered, unoriented Descartes configuration in the packing. Note that
$\bS_1^2=\bS_2^2=\bS_3^2=\bS_4^2=\bI$, 
and for any word $\bS_{i_1}\bS_{i_2}\dots \bS_{i_n}$,
($i_j \neq i_{j+1}$), the Descartes configuration
 $\sD'=(C_1', C_2', C_3', C_4')$ with ACC coordinates
$\bW_{\sD'} =\bS_{i_1}\bS_{i_2} \dots \bS_{i_n} \bW_\sD$
is a positively oriented Descartes configuration containing a unique circle
  $C_{i_n}'$ in the $n$-th stage partial packing.
Since no two circles at level $n$ touch for $n \geq 1$, by induction
it is easy to see that
\[
\sA=\{ \bS_{i_1} \bS_{i_2} \dots \bS_{i_n} \ | \ n \geq 0, i_j \neq i_{j+1} \},
\]
and
\[
\sA(\sD)=\{ \sD' \ | \  \bW_{\sD'}=\bS_{i_1} \bS_{i_2} \dots 
\bS_{i_n} \bW_\sD ,
\ n \geq 0,
i_j \neq i_{j+1} \}.
\]
Furthermore, for any ordered, positively oriented configuration $\sD''$ in
$\sP$, there is a unique permutation $\sigma$ such that $\sigma(\sD'')=\sD'
\in \sA(\sD)$, where there is a  unique sequence $i_1i_2\dots, i_n$,
$1 \leq i_j \leq 4$ and $i_j \neq i_{j+1}$, such that
$\bW_{\sD'} = \bS_{i_1} \bS_{i_2} \dots \bS_{i_n} \bW_{\sD}$.
This proves that all the ordered, positively oriented Descartes configurations
in an Apollonian circle packing comprise a union of 24 orbits of the Apollonian
group. The theorem follows by counting both the positively and negatively
oriented Descartes  configurations.
 $~~~\bsq$ \\

Theorem~\ref{th41}  shows  that the orbit under
the Apollonian group of a single (ordered, oriented)
Descartes configuration completely describes an Apollonian packing.

%
%
%
\section{Dual Apollonian Group}

As explained in \S3, the operation of 
inversion in the each of 
the individual circles of a Descartes configuration
is described by an integral matrix in $Aut(Q_{D})$.
This leads us to the following definition.

\begin{defi}\label{de43}
{\rm The {\em dual Apollonian group} $\sA^{\perp}$
is the subgroup of $Aut(Q_{D})$ generated by the
matrices 
\beql{S413}
\begin{array}{rllcrll}
\bS_1^\perp  & = & \left[
\begin{array}{rccc}
-1 & 0 & 0 & 0 \\
2 & 1 & 0 & 0 \\
2 & 0 & 1 & 0 \\
2 & 0 & 0 & 1
\end{array}
\right] ~, & \qquad & \bS_2^\perp & = & \left[
\begin{array}{crcc}
1 & 2 & 0 & 0 \\
0 & -1 & 0 & 0 \\
0 & 2 & 1 & 0 \\
0 & 2 & 0 & 1
\end{array}
\right] ~, \\ [+.3in]
\bS_3^\perp & = & \left[
\begin{array}{ccrc}
1 & 0 & 2 & 0 \\
0 & 1 & 2 & 0 \\
0 & 0 & -1 & 0 \\
0 & 0 & 2 & 1
\end{array}
\right] ~, && \bS_4^\perp & = & \left[
\begin{array}{cccr}
1 & 0 & 0 & 2 \\
0 & 1 & 0 & 2 \\
0 & 0 & 1 & 2 \\
0 & 0 & 0 & -1
\end{array}
\right] ~.
\end{array}
\eeq
}
We note that $\bS_i^{\perp}= \bS_i^{T}$, the transpose of $\bS_i$.
\end{defi}

By analogy with
Theorem~\ref{th41} we  might
think of  an orbit of the dual Apollonian group acting on
a single oriented Descartes configuration $\sD$,
as a  ``dual Apollonian circle packing''.
It can be viewed algebraically as a set of Descartes configurations
or geometrically as a collection of circles in the plane. 
In the algebraic viewpoint the orbit is a discrete set of points in
the parameter space $\sM_{\dd}$ so is a discrete object. 
From  the  geometrical viewpoint,
as a collection of circles, which we denote $\sP_{\sD}^{\perp}$,
it has the 
following weak ``packing''  property.
%
%

\begin{theorem}~\label{th41a} 
No  two circles in a dual Apollonian circle packing
cross each other. That is, two circles  in distinct
Descartes configurations in a dual Apollonian circle packing
$\sP_{\sD}^{\perp}$
either coincide, or are tangent to each other, or
are disjoint.
\end{theorem}

We do not give a proof of this theorem here, as it 
follows from a similar result proved for the
super-Apollonian packing in \S3 of part II.
It also is a special case of an $n$-dimensional generalization
proved  in \S4.2 of part III.

\begin{figure}[htbp]
\centerline{\epsfxsize=5.0in \epsfbox{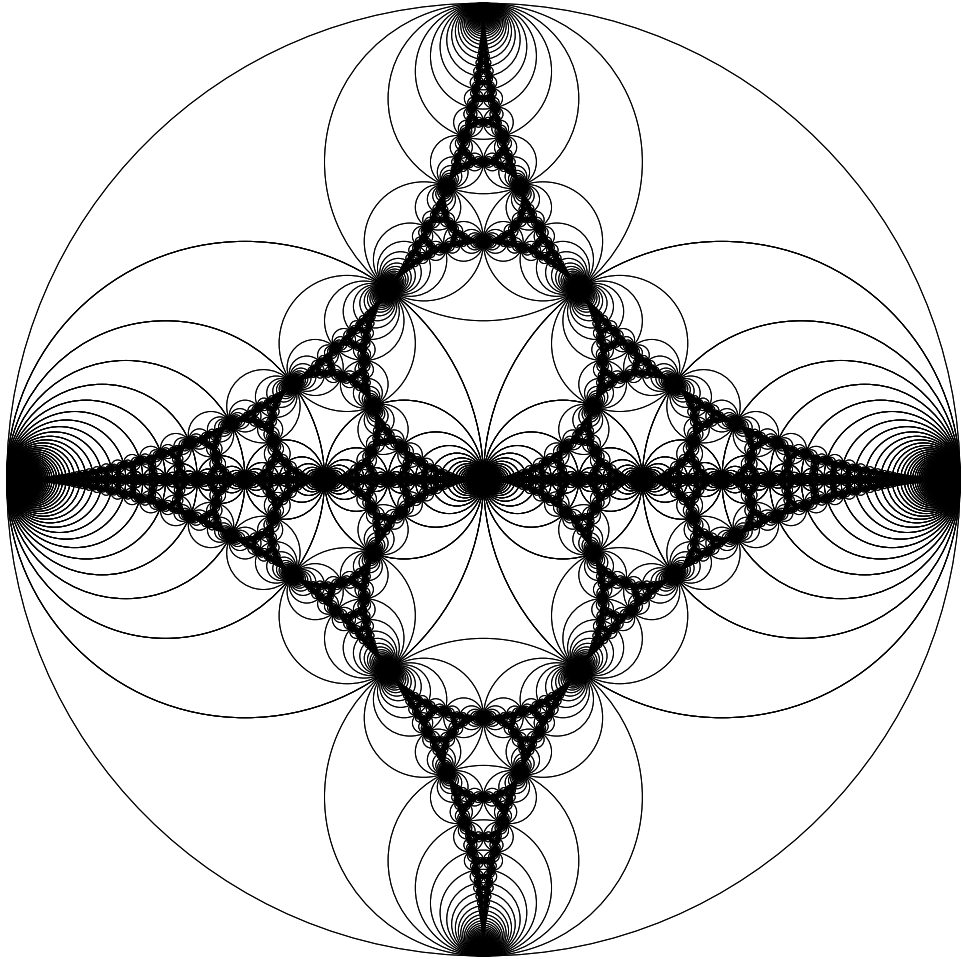}}
\caption{A dual Apollonian packing}\label{dualpack}
\end{figure}

Figure~\ref{dualpack}
pictures the  circles in a dual packing $\sP^{\perp}$. 
These are circles in Descartes configurations dual
to those in the Apollonian packing in Figure \ref{std}.
The fractal-like part of this figure is the closure of
the set of points at which circles  in the dual packing touch.
It can be shown that  infinitely many circles touch at each such tangency
point. Below we shall show that this fractal set coincides with the
limit set of the Apollonian packing generated by the
dual Descartes configuration $\bD(\sD_0)$. 
%
%

The  dual Apollonian group has a
simple relation  to the Apollonian group.

%
%

\begin{theorem}~\label{th42}
The dual Apollonian group $\sA^{\perp}$ is conjugate to the Apollonian group
$\sA$ using  the duality operator $\bD \in Aut(Q_D)$. This holds
at the level of generators, with 
$$
\bD^T \bS_i \bD = \bS_i^{\perp} ~~\mbox{for}~ 1 \le i \le 4.
$$
where $\bD= \bD^{-1} = \bD^T$.
\end{theorem}

\paragraph{Proof.} This is a straightforward computation.  $~~~\bsq$

%
%

\begin{theorem}~\label{th43}
Given a dual Apollonian packing 
$\sP_{\sD_0}^{\perp}$ generated by Descartes configuration $\sD_0$, 
let $\Sigma$
be the set of all intersection points
of the circles in $\sA^{\perp}(\sD_0)$. 
Then the  closure $\Lambda := \overline{\Sigma}$
 is equal to the residual set $ \Lambda(\sP_{\bD(\sD_0)})$ 
of the Apollonian packing $\sP_{\bD(\sD_0)}$ generated by
the dual Descartes configuration
$ \bD(\sD_0)$.
\end{theorem}

\paragraph{Proof.}
Theorem~\ref{th41b}(3) states
that   the residual set $\Lambda(\sP)$ of a
geometric  Apollonian
packing $\sP$ is the closure of the 
set of tangency points $\Sigma'$ in $\sP$.
Thus it will suffice to show that the set of
tangency points $\Sigma$ of the circles in 
the dual packing $\sP_{\sD_0}^{\perp}$
coincides with the set of tangency points
$\Sigma'$ of circles in the Apollonian packing
$\sP_{\bD (\sD_0)}$.
Each intersection point of two circles in the dual
Apollonian packing $\sP_{\sD_0}^{\perp}$
is an intersection point of circles in some Descartes configuration
$\sD$ of $\sA^{\perp}(\sD_0)$. (This follows from the recursive
construction of the dual packing.)
Each intersection point or circles in $\sD$ is 
an intersection point of circles in the dual Descartes 
configuration $\bD(\sD)$, which is the configuration
having augmented curvature-center coordinates
$\bD \bW_{\sD}$. Now $\bD(\sD)$
belongs to the Apollonian packing
$\sA(\bD (\sD_0))$ because
$$
\bD \bW_{\sD} = \bD (\bS_{i_1}^{\perp} \cdots \bS_{i_m}^{\perp} \bW_{\sD_0})=
\bS_{i_1} \cdots \bS_{i_m} (\bD \bW_{\sD_0}),
$$
using Theorem~\ref{th42}. Thus the set $\Sigma$ 
of intersection points
of circles  in the dual Apollonian packing 
$\sP_{\sD_0}^{\perp}$ is contained in 
the set $\Sigma'$ of intersection points of circles in 
the Apollonian packing $\sP_{\bD(\sD_0)}$.
The argument reverses to
show the converse is also true: each
intersection point of circles in the Apollonian
packing $\sP_{\bD (\sD_0)}$ is an intersection point of circles
in the dual Apollonian packing  $\sP_{\sD_0}^{\perp}$.
Thus $\Sigma = \Sigma'$. $~~~\bsq$

We have the following dichotomy between 
the geometric and algebraic views of Apollonian
packings and dual Apollonian packings.
Viewed geometrically as collections of circles,
a dual Apollonian packing $\sP^{\perp}$ 
appears quite  different from an Apollonian
packing, as the circles are nested to an infinite depth. 
However when viewed  algebraically as a collection of 
ordered, oriented Descartes
configurations $\DD(\sP^{\perp})$, it is a discrete set
given by the orbit of a group $\sA^{\perp}$ conjugate
to the Apollonian group $\sS$, and so is a similar object.

%
%
%

\section{Super-Apollonian group}
\setcounter{equation}{0}

The Apollonian group $\sA$ and dual Apollonian group $\sA^{\perp}$
both consist of integer automorphs in $Aut(Q_{D})$. We 
obtain a larger integral group by combining the
two groups, as follows.

\begin{defi}\label{de44}
{\rm The {\em super-Apollonian group} $\sA^{S}$ is the
subgroup of $Aut(Q_{D})$ generated by the Apollonian
group $\sA$ and $\sA^{\perp}$ together. We have
$$
\sA^{S} = \langle \bS_1, \bS_2, \bS_3, \bS_4, \bS_1^{\perp},  \bS_2^{\perp},
\bS_3^{\perp},  \bS_4^{\perp} \rangle.
$$
}
\end{defi}

The super-Apollonian group  $\sA^{S}$ is a discrete subgroup
of $Aut(Q_D)$ as a consequence of the fact that all its members
are integral matrices.

\subsection{Presentation of super-Apollonian group}

We now determine a presentation of
 the  super-Apollonian group, showing that it  
is a hyperbolic Coxeter group.
This result implies that both the 
Apollonian group and dual Apollonian group 
are finitely presented, and are hyperbolic 
Coxeter groups.

\begin{theorem}~\label{Sth47}
The  super-Apollonian group
$$
\sA^{S} = \langle \bS_1,~\bS_2,~\bS_3,~\bS_4,
~\bS_1^{\perp},~\bS_2^{\perp},~\bS_3^{\perp},~\bS_4^{\perp} \rangle
$$
is a hyperbolic Coxeter group
whose complete set of Coxeter relations are
\beql{S422b}
\bS_i^2 = ~(\bS_i^\perp)^2 = ~\bI \quad\mbox{for}\quad 1 \leq i \leq 4.
\eeq
\beql{S423}
(\bS_i \bS_j^\perp)^2~= ~(\bS_j^\perp \bS_i)^2 = ~\bI \quad\mbox{if}\quad i \neq j.
\eeq
\end{theorem}

\paragraph{Proof.}
Recall that $\bS_i^\perp = \bS_i^T.$
The group $\sA^{S}$ satisfies the relations
$\bS_i^2 = (\bS_i^T)^2 =\bI$ for $1 \le i \le 4$ and it satisfies
\beql{U429}
\bS_i \bS_j^T = \bS_j^T \bS_i \quad\mbox{for}\quad
i \neq j ~,
\eeq
which is equivalent to
\eqn{S423}.

Write words $\bU$ in $\sA^{S}$ as $\bU = \bU_n \bU_{n-1} \cdots \bU_1$,
 in which each $\bU_i$ is a generator $\bS_j$ or $\bS_j^T$, with the empty
word (n = 0) being the identity element.
We do not need inverses since each generator is its own inverse.
A word $\bU$ is in {\em normal form} if it has the two properties:
\begin{itemize}
\item[(i)]
$\bU_{k} \neq \bU_{k - 1} $ for $2 \le k \le n$.
\item[(ii)]
If $\bU_k = \bS_j^T$ for some $j$, then 
$\bU_{k - 1} \neq \bS_i$ for all $i \neq j$.
\end{itemize}
We can reduce any word to a word in normal form, using the relations
to move any symbols $\bS_j^T$ as far to the right in the word
as possible.  We cancel any adjacent identical symbols to make (i) hold.
Then move the rightmost $\bS_j^T$ as far to the right as possible,
as allowed by \eqn{U429}.
Repeat the same with the next rightmost symbol $\bS_j^T$.
If in the process any adjacent symbols $(\bS_j^T)^2$ occur, cancel them.
This process must terminate in a normal form word or the empty word.

The theorem is equivalent to showing that no normal form word with
$n \ge 1$ is the identity element in $\sA^{S}$.
This holds because the
 reduction to normal form used only the Coxeter relations,
so if $\sA^{S}$ satisfies an additional nontrivial relation,
there would exist some nontrivial normal form word that is the identity.

We proceed by induction on the length $n$ of a word in normal form
$\bU =\bU_n \bU_{n-1} \cdots \bU_1$.
Note that any suffix $\bU_{j-1} \bU_j \cdots \bU_1$ 
of $\bU$ is also in normal form.
For each $n \ge 1$, let $\sX_n$ denote
the set of  all  normal form words of length $n$.
(By a simple enumeration we can show that  the number of normal form words 
of length $n$ is $9 \cdot 5^{n-1}-1$. ) 
We measure the {\em size} of a word $\bU$ 
(viewed as a $4 \times 4$ matrix) as
\beql{U430}
f(\bU) := \bo^T \bU \bo = \sum_{i=1}^4 \sum_{j=1}^4 \bU_{ij} ~.
\eeq
Thus $f(\bI) = 4.$ For $n=1$ we have
\beql{U431}
f(\bU) = f(\bU_1) = 8
\eeq
in all cases.
For $n \ge 2$ if $\bU \in \sX_n$ then
$\bU' = \bU_{n-1} \bU_{n-2} \cdots \bU_1 \in \sX_{n-1}$, and we will  prove
\beql{U432}
f(\bU) > f(\bU' ) ~.
\eeq
If so, then
$$f(\bU) \ge f(\bU_1) \ge 8, $$
hence $\bU \neq \bI$, which will complete the proof.

We let
\beql{U433}
\br (\bU) := (r_1 ( \bU) , r_2 (\bU) , r_3 (\bU) , r_4 (\bU))^T = \bU \bo \,,
\eeq
be the vector of row sums of $\bU$. Now
\begin{eqnarray*}
f(\bS_1 \bU') & = & \bo^T
\bS_1 \br (\bU') = (-1,~ 3,~ 3,~ 3) \br (\bU') = -4 r_1 (\bU') + 3 f(\bU' ) \,, \\
f(\bS_1^T \bU') & = & \bo^T \bS_1^T \br (\bU') = (5,~1,~1,~1) \br (\bU)
= 4r_1 (\bU') + f(\bU'),
\end{eqnarray*}
with similar formulas in the other cases.
To prove \eqn{U432} it is therefore sufficient to prove the
following two assertions,
for $1 \le h \le 4$.
\begin{itemize}
\item[(1)]
If $\bU_n = \bS_h$, so that $\bU_{n-1} \neq \bS_h$, then 
$2r_h (\bU') < f(\bU')$.
\item[(2)]
If $\bU_n = \bS_h^T$, so that $U_{n-1} \neq \bS_h^T$ and 
$U_{n - 1} \neq \bS_i$
for $i \neq h$,
then $r_h (\bU') > 0$.
\end{itemize}

\noindent

Instead of proving (1), (2), we
 will prove by induction on $n$ the following three assertions for all
$\bU \in \sX_n$. Here  $(h,i,j,k)$ always denotes some permutation of
$(1,2,3,4)$ in what follows.

\begin{itemize}
\item[(i)]
For $1 \le i,j \le 4$, with $i \neq j$,
$$r_i (\bU) + r_j (\bU) > 0 ~.$$
\item[(ii)]
If $\bU_n = \bS_h$, then
\begin{eqnarray*}
r_h (\bU) & > & 0, \\
r_i (\bU) & < & r_h (\bU) + r_j (\bU) + r_k (\bU), \\
r_j (\bU) & < & r_h (\bU) + r_i (\bU) + r_k (\bU), \\
r_k (\bU) & < & r_h (\bU) + r_i (\bU) + r_j (\bU)~.
\end{eqnarray*}

\item[(iii)]
If $\bU_n = \bS_h^T$ then
\begin{eqnarray*}
r_h (\bU) & < & 0, \\
r_i (\bU) & < & r_h (\bU) + r_j (\bU) + r_k (\bU), \\
r_j (\bU) & < & r_h (\bU) + r_i (\bU) + r_k (\bU), \\
r_k (\bU) & < & r_h (\bU) + r_i (\bU) + r_j (\bU) ~.
\end{eqnarray*}
\end{itemize}

\noindent
Note that (i) implies that at most one of the row sums of $\bU$ 
can be negative.
If proved, (i)--(iii) imply (1), (2), which themselves imply \eqn{U432},
completing the proof of the Theorem.

The induction hypothesis (i)--(iii) holds for $n=1$, since
$r(\bS_h)$ is a permutation of $(5,1,1,1)^T$ and $r(\bS_h^T)$
is a permutation of $(-1,3,3,3)^T$.
Suppose it is true for $n$.
For the induction step, write $\bU \in \sX_{n+1}$ as
$\bU = \bU_{n+1} \bU'$ with $\bU' \in \sX_n$ 
and abbreviate $r'_i = r_i (\bU')$.

\paragraph{Case 1.}
{\em $U_{n+1} = \bS_h$ for some $h$.}

\noindent
We must verify (i) and (ii) for $\bU$.
We have
$$
\begin{array}{lllll@{}l}
r_h (\bU) & = & r_h (\bS_h \bU') & = & - & r'_h + 2r'_i + 2r'_j + 2r'_k, \\
r_i (\bU) & = & r_i (\bS_h \bU') & = & r'_i, \\
r_j (\bU) & = & r_j (\bS_h \bU') & = & r'_j, \\
r_k (\bU) & = & r_k (\bS_h \bU') & = & r'_k ~.
\end{array}
$$
To verify (i), all cases not involving $r_h (\bU)$ follow from the induction
hypothesis.
To show $r_h (\bU) + r_i (\bU) > 0$ note that
\beql{U433a}
r_h (\bU) + r_i (\bU) = - r'_h + 3r'_i + 2r'_j + 2r'_k
~.
\eeq
Now $\bU_{n} \neq \bS_h$ since $\bU \in \sX_{n+1}$, and for all cases
except $\bU_{n} = \bS_h^T$ the induction hypotheses (ii) and (iii) give
$r'_h < r'_i + r'_j + r'_k$, whence
$$r_h (\bU) + r_i (\bU) \ge 2r'_i + r'_j + r'_k =
(r'_i + r'_j ) + (r'_i + r'_k) > 0 \,,
$$
using induction hypothesis (i).
If $\bU_{n} = \bS_h^T$, then by hypothesis (iii) $r'_h < 0$,
hence $r'_i, r'_j , r'_k > 0$ by
hypothesis (i), so all terms on the right side of \eqn{U433a} are positive,
so $r_h (\bU) + r_i (\bU) > 0$.
The cases
$r_h (\bU) + r_j (\bU) > 0$ and $r_h (\bU) + r_k (\bU) > 0$ follow similarly.
To verify (ii) for $\bU$ it suffices to prove
\begin{eqnarray*}
&& 0 < - r'_h + 2r'_i + 2r'_j + 2r'_k, \\
&& r'_i < -r'_h + 2r'_i + 3r'_j + 3r'_k, \\
&& r'_j < - r'_h + 3r'_i + 2r'_j + 3r'_k, \\
&& r'_k < - r'_h + 3r'_i + 3r'_j + 2r'_k,
\end{eqnarray*}
which is equivalent to
\begin{eqnarray}\label{U434}
r'_h & < & 2r'_i + 2r'_j + 2r'_k, \nonumber \\
r'_h & < & r'_i + 3r'_j + 3r'_k, \nonumber \\
r'_h & < & 3r'_i + r'_j + 3r'_k,  \\
r'_k & < & 3r'_i + 3r_j + r'_k ~. \nonumber
\end{eqnarray}
There are three cases, according as
$\bU_n = \bS_i$, $\bU_n = \bS_i^T$ with $i \neq h$ and $\bU_n = \bS_h^T$.
In the first two of these,
the induction hypotheses give
$$r'_h < r'_i + r'_j + r'_k,$$
which with induction hypothesis (i) immediately yields \eqn{U434}.
If $\bU_n = \bS_h^T$ then by (iii), $r'_h < 0$, whence
$r'_i , r'_j , r'_k > 0$ by (i) so \eqn{U434} is immediate.
This finishes Case 1.

\paragraph{Case 2.}
{\em $\bU_{n+1} = \bS_h^T$ for some $h$.}

\noindent
We must verify (i) and (iii).
We have
$$
\begin{array}{lllll}
r_h (\bU) & = & r_h (\bS_h^T \bU' ) & = & -r'_h, \\
r_i (\bU) & = & r_i (\bS_h^T \bU') & = & 2r'_h + r'_i, \\
r_j (\bU) & = & r_j (\bS_h^T \bU') & = & 2r'_h + r'_j, \\
r_k (\bU) & = & r_k (\bS_h^T \bU') & = & 2r'_h + r'_j,
\end{array}
$$
To prove (i), note first that
$$r_h (\bU) + r_i (\bU) = r'_h + r'_i > 0,$$
using (i) for $\bU'$, and similar inequalities hold
 for other cases involving $r_h (\bU)$.
The remaining cases are all of the form
$r_i (\bU) + r_j (\bU) > 0$
with $i,j \neq h$, and the proofs for each are similar.
We have
$$
r_i (\bU) + r_j (\bU) = 4r'_h + r'_i + r'_j =
2r'_h + (r'_h + r'_i) + (r'_h + r'_j ) > 2r'_h \,,
$$
using induction hypothesis (i).
We now show $r'_h > 0$ holds in all cases.
Since $\bU \in \sX_{n+1}$, we have
$\bU_{n} \neq \bS_h^T$ and $\bU_{n} \neq \bS_i$ with $i \neq h$.
If $\bU_{n} = \bS_j^T$ for some $j \neq h$ then hypothesis (iii) for
$\bU'$ says that $r'_j < 0$, whence we must have $r'_h > 0$.
If $\bU_{n} = \bS_h$ then hypothesis (ii) for $\bU'$ says $r'_h > 0$.
Thus $r'_h > 0$ in all cases, and (i) holds for $\bU$.
To prove (ii) we must prove
$$
\begin{array}{l}
r'_h < 0, \\
2r'_h + r'_i < 3r'_h + r'_j + r'_k, \\
2r'_h + r'_j < 3r'_h + r'_i + r'_k, \\
2r'_h + r'_k < 3r'_h + r'_i + r'_j,
\end{array}
$$
which is equivalent to
\beql{U435}
\begin{array}{l}
r'_h > 0, \\
r'_i < r'_h + r'_j + r'_k, \\
r'_j < r'_h + r'_i + r'_k, \\
r'_k < r'_h + r'_i + r'_j ~.
\end{array}
\eeq
We verified $r'_h > 0$ already.
For the remainder there are two cases, according as
$\bU_n = \bS_h$ or $\bU_n = \bS_i^T$ for some $i \neq h$.
If $\bU_n = \bS_h$ the three remaining inequalities in \eqn{U435}
follow from inductive hypothesis (ii) for $\bU'$.
If $\bU_n = \bS_i^T$ the inductive hypothesis (iii) gives the last two
inequalities in \eqn{U435}, and also that $r'_i < 0$.
By hypothesis (i) for $\bU'$ this implies $r'_j , r'_h , r'_k > 0$ hence
$r'_i < 0 < r'_h + r'_j + r'_k$ which verifies \eqn{U435} in this case,
and finishes Case 2.

The induction is complete.~~~$\bsq$
%
%
%

\subsection{Apollonian super-packings}

We now define an {\em Apollonian super-packing}
to be a set of (ordered, oriented) Descartes configurations forming   
an orbit of the super-Apollonian group $\sA^{S}$,
acting on a single such configuration. 

The study of Apollonian super-packings will form the subject matter
of part II. We call the set of circles in all the
Descartes configurations in such a super-packing
a {\em geometric Apollonian super-packing}. 
These circles comprise four orbits of a group of
M\"{o}bius transformations
$$
\sG_{\sA^{S}}(\sD) = \langle \fs_1, \fs_2, \fs_3, \bs_4, \bs_1^{\perp},
\bs_2^{\perp}, \bs_3^{\perp}, bs_4^{\perp}\rangle 
$$
where $\fs_i$ and $\fs_i^{\perp}$
were defined in \S3.3. In terms of the space
$\sM_{\DD}$ of Descartes configurations this group action
is given by 
$\tilde{\sG}_{\sA^{S}}(\sD) \subset Aut(Q_{W})$ 
defined by 
$$
\tilde{\sG}_{\sA^{S}}(\sD) := \bW_{\sD}^{-1}  \sA^{S} \bW_{\sD},
$$
using the isomorphism $Aut(Q_{W}) = \mbox{M\"{o}b}(2) \times \{ \bI, -\bI\}$
given in Theorem~\ref{th91} of Appendix A.

In part II we  will show that 
the individual circles in a geometric super-packing
form a ``packing'' in
the weak sense that no two circles cross each
other, although circles can be nested.
This is a remarkable geometric fact, because 
we will also show in part II that a 
strongly integral Apollonian super-packing necessarily contains a copy of
every integral Apollonian circle
packing, with all these copies (essentially)
 contained inside the square $0 \le x \le 2, 0 \le y \le 2$.
All these Apollonian packings
fit together in such a way that all the circles in these packings 
manage not to cross each other.

%
%
%
\newpage
\section{Appendix A.  M\"{o}bius Group Action}
\setcounter{equation}{0}
The {\em (general) M\"{o}bius group} 
 $\mbox{M\"{o}b}(2)$
is the group of M\"{o}bius transformations, 
allowing reflections
(including complex conjugation).  This group is denoted $GM(\hat{\rr}^2)$ in 
Beardon~\cite[p. 23]{Be83},
and is also known as the {\em conformal group}. 
The group
$\mbox{M\"{o}b}(2)_{+}  \simeq PSL(2, \cc) = SL(2, \CC)/\{\pm \bI \}$
consists of the orientation-preserving maps of 
$\hat{\cc}= \rr^2 \cup \{\infty\}$.
The group $\mbox{M\"{o}b}(2)$
is  a six-dimensional real Lie group which has two connected
components corresponding to orientation-preserving
\footnote{The terminology ``orientation'' here 
refers to the invariant $\det(\bM)$.} and
orientation-reversing transformations.
It can be written as a semi-direct product
$$
\mbox{M\"{o}b}(2) \cong \mbox{M\"{o}b}(2)_{+} \rtimes \{ 1, \fc\},
$$
in which $\fc$ denotes complex conjugation.
M\"{o}bius transformations
take circles to circles (or straight lines)
and preserve angles.
Thus they take ordered Descartes configurations to ordered 
Descartes configurations. Concerning orientation, we show
that although  M\"obius transformations can reverse orientation of a
single circle, they preserve (total) orientation of 
oriented Descartes configurations.
Recall that the {\em total orientation} 
of an oriented Descartes configuration is the sign of the sum of the
(signed) curvatures of the circles in it.

\begin{theorem}~\label{th90a}
M\"obius transformations preserve the total orientation of
oriented Descartes configurations.
\end{theorem}

\pf
Positively oriented Descartes configurations are characterized
by the four circles having disjoint ``interiors'', specified  by
the positive orientation. 
M\"obius transformations take Descartes configurations
to Descartes configurations, and preserve the ``disjoint
interior'' property, hence preserve positive orientation, as
given by normal vectors to the circle. The result for
negatively oriented Descartes configurations holds since
all normal vectors are reversed from the positively oriented  case.
~~~$\bsq$ \\

We now consider a group $GM^{*}(2)$
with four connected components,
which we will term 
the  {\em extended General M\"{o}bius group},
defined by 
$GM^{*}(2) := \mbox{M\"{o}b}(2) \times \{ \bI, - \bI\}$.
Here $\{ \bI, -\bI\}$ are in the center of this group,
and we  write elements of $GM^{*}(2)$ as $\pm \fg$,
in which  $\fg \in \mbox{M\"{o}b}(2)$, 
and the sign indicates
which of  $\pm \bI$ occurs.
We have 
$$
GM^{*}(2) \cong SL(2, \CC) \rtimes \{ \bI, \fc \}.
$$

%
%

\begin{table}[t]
 \begin{tabular}[t]{|c|c|}
 \hline
  \begin{diagram}  
   GM^{*}(2)=  M\ddot{o}b(2)\times \{-\bI, \bI\}    & \rTo^\sim & Aut(Q_W)  & 
         \rTo^\sim  &   {O(3,1)}  \\
   \dTo^\pi \uInto &   &  \dTo^\pi \uInto &    & \dTo^\pi \uInto  \\
   M\ddot{o}b(2)=PSL(2, \CC)\rtimes \{I, \fc\} & \rTo^\sim & 
Aut(Q_W)^\uparrow 
         &  \rTo^\sim & O(3,1)^\uparrow  \\ 
   \uInto &  &  \uInto& & \uInto  \\
   {M\ddot{o}b(2)_+=PSL(2,\CC) }      & \rTo^{\sim} 
         & Aut(Q_W)^\uparrow_+ &  \rTo^\sim  &  O(3,1)^\uparrow_+  \\ 
  \end{diagram}
 &  
  \begin{diagram} 
  \text{Lorentz group} \\ 
   {} \\
  {\text{orthochronous} \atop  \text{Lorentz group}} \\ 
   {}\\
  {\text{proper orthochronous}  \atop \text{Lorentz  group}} \\
 \end{diagram} 
 \\
 \hline
\end{tabular}
\caption{Group Isomorphisms}~\label{table1}
\end{table}

The main object of this Appendix is to  define 
an action of $GM^{*}(2)$ on the right on the parameter space $\sM_{\dd}$,
given in the next theorem. This amounts to
finding an explicit isomorphism between $GM^{*}(2)$
and $Aut(Q_{W})$, which appears as the  horizontal
arrow  on the left in the top row in Table ~\ref{table1}.
This map when restricted to the smaller groups
M\"{o}b(2) and M\"{o}b$(2)_{+}$ give the other
two horizontal isomorphisms on the left side of
the table. 
Table A-1 also indicates isomorphisms on its right side  
to the
Lorentz group $O(3,1)$ and 
corresponding  subgroups, which we defer discussing 
until  after the following result.

%
%
%

\begin{theorem}\label{th91}
Let $GM^{*}(2) := \mbox{M\"{o}b}(2) \times \{ \bI, - \bI\}$.
There is a unique 
isomorphism $\pi: GM^{*}(2) \to Aut(Q_W)$,
with image elements $ \bV_{\pm \fg}:=\pi(\pm \fg)$,
such that the following hold.
 
(i) For  $\fg \in  \mbox{ M\"{o}b}(2)$
the augmented curvature-center coordinates
for each ordered, oriented Descartes configuration $\sD$ 
satisfy 
\beql{N9.2}
\bW_{\fg(\sD)} = \bW_{\sD} \bV_{\fg}^{-1}.
\eeq

(ii) The action of $-\bI$ on augmented
curvature-center coordinates is
\beql{N9.2b}
\bW_{-\sD} = \bW_{\sD}\bV_{-\bI}^{-1}= -\bW_{\sD}.
\eeq 
\end{theorem}

\paragraph{Proof.}
We compute the action of 
$\mbox{M\"{o}b}(2)$
acting on
augmented curvature-center coordinates.
Let $(\bar{b}, b, w_1, w_2)= (\frac{x_1^2+y_1^2-r^2}{r}, 
\frac{1}{r}, \frac{x_1}{r}, \frac{y_1}{r})$
be the augmented curvature-center coordinates of the circle
$$ (x - x_1)^2 + (y - y_1)^2 = r^2. $$
This circle can be recovered from these coordinates via
\beql{281}
(bx - w_1)^2 + (by - w_2)^2 = 1,
\eeq
and the orientation of the circle (inside versus outside) is
determined by the sign of $b$. An oriented ``circle at infinity''
is a line given by 
\beql{282}
x\cos{\theta}  +y \sin{\theta}  = m,
\eeq
and its associated curvature-center coordinates are 
\beql{283}
(\bar{b}, b, w_1, w_2)= (2m, 0,\cos{\theta}, \sin{\theta}).
\eeq
 Here the
orientation is given by the convention that the  
normal $(\cos{\theta},\sin{\theta})$ points inward.

The group M\"{o}b(2) is generated by 

 (1) translations $\ft_{z_0}= z + z_0$;

 (2) dilations $\fd_\lambda(z) = \lambda z$ with $\lambda \in \cc$, 
$\lambda \neq 0$;  

 (3) the conjugation $\fc (z) = \bar{z}$;

 (4)  the inversion in the unit circle
$\fj(z) = \frac{1}{\bar{z}} = \frac{z}{|z|^2}.$

Given $\fg \in \mbox{M\"{o}b}(2)$,
we will let $\tilde{\fg}$ denote the 
corresponding action on  the curvature-center coordinates of an oriented
circle. 
The action of translation by $z_0 = x_0 + iy_0$ is
\beql{284}
\tilde{\ft}_{z_0}(\bar{b}, b, w_1, w_2) = 
(\bar{b}+2w_1x_0+2w_2y_0+b(x_0^2+y_0^2), b, w_1 + bx_0, w_2 + by_0).
\eeq
The action of a dilation with $\lambda=r e^{i\theta}$, $(r >0)$ is given by 
\beql{285}
\tilde{\fd}_\lambda (\bar{b}, b, w_1, w_2) = 
(r\bar{b}, \frac{b}{r},~ w_1\cos{\theta}-w_2\sin{\theta}, w_1\sin{\theta}+w_2\cos{\theta}).
\eeq
The action  of complex conjugation is
\beql{286}
\tilde{\fc}(\bar{b}, b, w_1, w_2) = (\bar{b}, b, w_1, - w_2).
\eeq
The action of inversion is
$$
\tilde{\fj}( \bar{b}, b, w_1, w_2) = (b, \bar{b}, w_1, w_2).
$$
All of these actions apply
to ``circles at infinity'' and extend  to linear maps on the $4 \times 4$ 
matrices $\bW_{\sD}$. 

The translation operation is given by right multiplication by
the matrix 
\beql{Gt0}
\bV_{\ft_{z_0}}^{-1} := \left[
\begin{array}{crrc}
1 &  0  &  0  &  0 \\
x_0^2+y_0^2  &  1   &  x_0   & y_0 \\
2x_0 &  0   &  1   & 0 \\
2y_0 & 0 & 0 & 1 
\end{array}
\right],
\eeq
and one verifies  \eqref{N9.2} holds by direct computation.

For the  dilation $\fd_\lambda$, with $\lambda = re^{i \theta}$
$(r > 0)$  the right action is by the matrix
\beql{Gdlambda}
\bV_{\fd_\lambda}^{-1} :=\left[
\begin{array}{cccc}
r         &   0            &     0          & 0 \\
0         & 1/r            &     0          & 0 \\
0         & 0              &  \cos{\theta}  & \sin{\theta}  \\
0         & 0              & -\sin{\theta}  & \cos{\theta}  
\end{array} 
\right].
\eeq
For complex conjugation
$\fc$, the  right action is by the matrix 
\begin{equation}
\bV_{\fc}^{-1} = \bV_{\fc} 
:=\left[
\begin{array}{rrrr}
1         &  0             &     0          & 0 \\
0         &  1             &     0          & 0 \\
0         &  0              &    1          & 0  \\
0         &  0              &    0          & -1 
\end{array} 
\right].
\end{equation}
For  the inversion $\fj_C$ in the unit circle, the permutation matrix 
\begin{equation}
\bV_{\fj_C}^{-1} =\bV_{\fj_C}= P_{(12)} = \left[
\begin{array}{cccc}
0         &   1            &     0          &    0 \\
1         &   0            &     0          &    0 \\
0         &   0            &     1          &    0  \\
0         &   0            &     0          &    1  
\end{array} 
\right].
\end{equation}

It is easy to verify that the above matrices are all in $Aut(Q_W)^{\uparrow}$, 
so that the map  so far defines  a homomorphism of $\mbox{M\"{o}b}(2)$
into $Aut(Q_W)^{\uparrow}\simeq O(3, 1)^{\uparrow}$,
identified with the  isochronous Lorentz group.
The group $\mbox{M\"{o}b}(2)$ acts simply transitively on
ordered Descartes configurations, as observed by Wilker
\cite[Theorem 3, p. 394]{Wi81},
and the group $Aut(Q_{W})$ acts simply transitively on
ordered, oriented Descartes configurations by Theorem~\ref{th33}.
Because  $Aut(Q_W)^{\uparrow}$ is of index $2$ in $Aut(Q_{W}) \simeq O(3,1)$,
we conclude that the map so far defines  an isomorphism 
of $\mbox{M\"{o}b}(2)$  onto $Aut(Q_W)^{\uparrow}$.

To complete the proof, we define  the action of $-\bI$ to be 
\begin{equation}~\label{flip}
(\bV_{-\bI})^{-1} = \bV_{-\bI} = - \bI.
\end{equation}
It has the effect of reversing  (total) orientation of the 
Descartes configuration, and does not correspond to 
a conformal transformation. Since $-\bI \notin Aut(Q_W)^{\uparrow}$,
adding it gives the desired isomorphism of
$GM(2)$ onto $Aut(Q_{W})$. $~~~\bsq$ \\

In   terms of  the natural coordinates
on $\mbox{M\"{o}b}(2)_{+} \simeq PSL(2, \CC)$ the homomorphism
$\pi(\cdot)$ given in Theorem~\ref{th91} 
is a non-linear map. This can be clearly seen
in \eqn{Gdlambda}, where both $\bU_{\fd_\lambda}$ and
$\bU_{\fd_\lambda}^{-1}$ are nonlinear functions of the
coordinates of
\beql{N817a}
{\fd}_{\lambda} = \left[ 
\begin{array}{cc}
\sqrt{\lambda} & 0 \\
0  & \frac{1}{\sqrt{\lambda}} 
\end{array} \right] \in SL(2, \CC).
\eeq
Here the  two choices of  $\pm \sqrt{\lambda}$ give a well defined
matrix modulo $\pm I$.

We now return to the data in  Table ~\ref{table1} 
giving  the isomorphisms of $Aut(Q_{W})$ and its
subgroups to the Lorentz group $O(3,1)$ and its
two subgroups 
$O(3,1)^{\uparrow}$ the
orthochronous Lorentz group, and $O(3,1)_{+}^{\uparrow}$
the proper orthochronous Lorentz group, 
using  the terminology
of Wilker ~\cite{Wi81}. An {\em orthochronous Lorentz transformation}
is one that preserves the arrow of time; Table ~\ref{table1} shows
this corresponds to preserving total orientation of a 
Descartes configuration.
The set of isomorphisms
given by the three horizontal arrows
on the right in Table~\ref{table1}
are obtained by any fixed  choice of matrix $\bA$
that intertwines $Q_{W}$ and $Q_{\sL}$ by
$\bQ_{W} = \bA^T \bQ_{\sL} \bA$, in which case the 
isomorphism is $Aut(Q_{W}) = \bA^{-1} O(3,1) \bA$
sending $\bV \mapsto \bA \bV \bA^{-1}$. 
To preserve the underlying rational structure
of these forms, one must choose $\bA$ to be a
rational matrix, for example the integral matrix
given in \eqn{N333b}.

There  is a different  choice of intertwining
matrix, for the rightmost horizontal rows in
Table A-1, which will allow us  to make the composed horizontal
maps in Table \ref{table1} from the M\"{o}bius
group to the Lorentz groups compatible with
the framework used by Wilker \cite{Wi81}.
In  \cite[Theorem 10]{Wi81} Wilker
gave an 
explicit isomorphism  of $PSL(2, \CC)$ onto 
the proper orthochronous Lorentz group $O(3,1)_{+}^{\uparrow}$
It maps
\footnote{We have permuted the first and last row and column
of Wilker's result because he uses the Lorentz form
$w^2 + x^2+y^2-z^2$.}
$$
\pm \left[ \begin{array}{rr}
a & b \\
c & d 
\end{array}
\right] 
\longmapsto
 \left[
\begin{array}{rrrr}
\frac{1}{2}\left(|a|^2 + |b|^2 \right. & \Im(a\bar{c} + b \bar{d}) & 
\frac{1}{2}\left (|a|^2 + |b|^2 \right. &  \Re(a\bar{c} + b\bar{d}) \\
\left. +|c|^2 + |d|^2\right ) & &  \left. - |c|^2-|d|^2 \right)& \\
\Im(-a\bar{b} - c\bar{d}) & \Re(a\bar{d} - b \bar{c}) &  \Im(-a\bar{b} +
c\bar{d}) &  \Im(-a\bar{d} + b \bar{c}) \\
\frac{1}{2}\left(|a|^2 - |b|^2 \right. &  \Im(a\bar{c}-b\bar{d}) & 
\frac{1}{2}\left(|a|^2 -|b|^2 \right. &  \Re(a\bar{c}-b \bar{d}) \\
\left. +|c|^2 - |d|^2\right)& & \left. - |c|^2 + |d|^2\right) & \\
\Re(a\bar{b}+c\bar{d}) & \Im(a\bar{d}+ b \bar{c}) &  \Re(a\bar{b}-c\bar{d}) &
 \Re(a\bar{d} + b \bar{c}) 
\end{array}
\right].
$$
One can identify this map  with the composed map in 
the last row of Table A-1 
from the first column  to the last column, if one 
chooses  the map $Aut(Q_{W}) \to Aut(Q_{\sL})$
to be  $\bV \mapsto \bZ \bV \bZ^{-1}$
for  a particular intertwining matrix $\bZ$,
which  is unique up to multiplication
by $\pm \bI$. It has irrational entries, and is 
\beql{convert}
\bZ := \pm \sqrt{2} \left[
\begin{array}{cccc}
1         &   1            &     0          &    0 \\
0         &   0            &     0          &    -1 \\
1         &   -1           &     0          &    0  \\
0         &   0            &     -1         &   0  
\end{array} 
\right].
\eeq

%
%
%
%
%
\newpage
\section{Appendix B.  Groups in Hyperbolic $3$-Space associated 
to Apollonian Packings} 
\setcounter{equation}{0}

In this appendix we  describe two different discrete
groups  of isometries of hyperbolic $3$-space $\HH^3$
associated to an Apollonian packing. We specify these
groups as groups of M\"{o}bius transformations of $\hat{\cc}$.

Any M\"{o}bius transformation of $\hat{\cc} = \hat{\rr}^2$
has a unique Poincare lift to an isometry of hyperbolic $3$-space
$\HH^3$, viewed as the upper half-space of 
$\hat{\rr}^3 = \rr^3 \cup \{\infty\}$,
as described in Beardon~\cite[Sect. 3.3]{Be83}. Here 
$\hat{\cc}$ is identified with the ideal boundary of $\hat{\rr}^3$.
In this way any  group of M\"{o}bius transformations acting on
$\hat{\cc}$  lifts to a group of hyperbolic isometries.

A {\em classical Schottky group} is a group 
$\Gamma = \langle \fg_1, ...\fg_n \rangle$ of  
M\"{o}bius transformations 
associated to $2n$ circles in the plane $\hat{\cc}$ having
disjoint interiors, but possibly touching on their boundaries.
The group has $n$ generators $\fg_j$,   in which the 
$j$-th generator of the group maps the exterior of circle $C_{2j-1}$
onto the interior of circle $C_{2j}$.
Schottky groups consist entirely of
holomorphic M\"{o}bius transformations.
(Note: knowing a single $\fg_j$ alone does not determine the
circles $C_{2j-1}$ and $C_{2j}$ uniquely.)

Mumford, Series and Wright \cite[Chapter 7]{MSW02}
observe that given an  Apollonian circle packing $\sP_{\sD}$
generated by a Descartes configuration $\sD$,
there is a Schottky group $\Gamma (\sD)$ in 
hyperbolic $3$-space $\HH^3$ whose limit
set is exactly the limit set $\Lambda_{\sD}$ 
of the Apollonian packing $\sP_{\sD}$.
The Schottky group $\Gamma(\sD) $ has two generators
$\bP_1, \bP_2$, involving maps between
two pairs of circles in the dual Descartes
configuration
$\sD^{\perp}$, specified by
the ordering as $A, a, B, b$,
say. It takes the exterior of circle $A$ to the
interior of circle $a$ and the exterior of circle $B$ to
the interior of circle $b$.  

\begin{figure}[htbp]
\centerline{\epsfxsize=4.0in \epsfbox{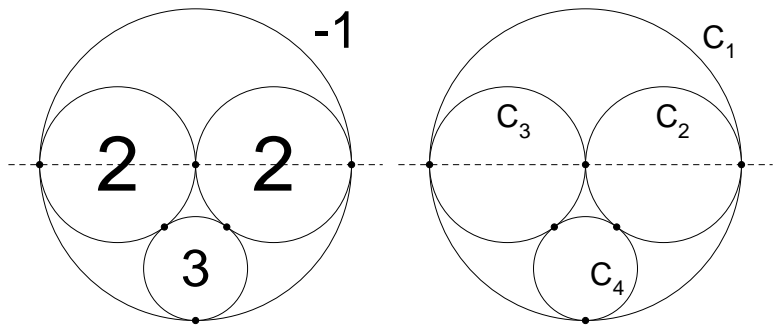}}
\caption{Special Descartes configuration $\sD_0$}~\label{B1-special}
\end{figure}

We treat the case they picture \cite[Figure 7.3]{MSW02}.
 where the Schottky limit
set is $\Lambda(\sD_0)$,
for  the (ordered) Descartes
configuration $\sD_0$ having coordinates
$$
\bW_{\sD_0}= \left[
\begin{array}{rrrr}
1  & -1 & 0 & 0 \\
0  & 2 & 1 & 0 \\
0  & 2 & -1 & 0 \\
1  & 3 &  0 & -2
\end{array}\right].
$$
This configuration is  pictured in Figure \ref{B1-special}, with part (a)
indicating curvatures and part (b) shows the ordering of  circles,
used later in labeling  M\"{o}bius  group generators $\fs_i$.

%
%
%

\begin{figure}[htbp]
\centerline{\epsfxsize=4.0in \epsfbox{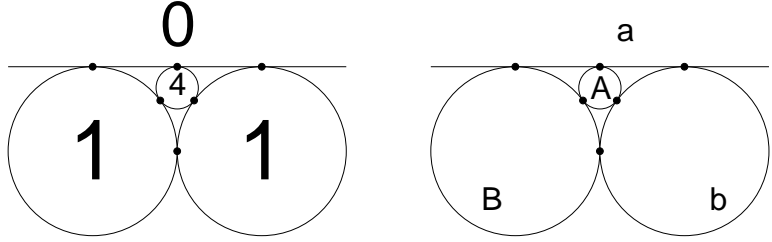}}
\caption{Dual Descartes configuration $\sD_0^{\perp}$}~\label{B2-dual}
\end{figure}

The Schottky group operations are associated to the dual
Descartes configuration $\sD_0^{\perp}$, which  has 
augmented curvature-center coordinates
$$
\bW_{{\sD}_0^{\perp}}= \left[
\begin{array}{rrrr}
0 & 0 & 0 & 1 \\
1 & 1 & -1 & -1 \\
1 & 1 & 1 & 1 \\
0 & 4 & 0 & -1 \\
\end{array}\right].
$$
It is pictured in Figure \ref{B2-dual}, with (a) indicating
curvatures and (b) the labeling of circles used in 
the Schottky group.
%
%
%

In Figure \ref{B2-dual}  the Schottky group generator 
$\bP_1$ maps the exterior of circle $A$ to the
interior of circle $a$, and $\bP_2$ maps the exterior
of circle $B$ to the interior of circle $b$.
We let $\Gamma_{S}(\sD)$ denote the Schottky group
they generate.
The limit set  $\Lambda_{\Gamma_{S}(\sD)}$ of
this Schottky group is the
closed set in  $\PP^1(\CC)$ obtained by removing 
the interiors of all circles in the Apollonian packing
associated to $\sD$. This is exactly
the residual set of the Apollonian packing,
 see Theorem~\ref{th41b} for its properties. It is
a fractal set of Hausdorff dimension about $1.3$.

The precise result is as follows.

\begin{prop}~\label{thB1}
The Schottky group $\Gamma = \langle \bP_1, \bP_2 \rangle$ given 
by the parabolic generators
\beql{1001}
\bP_1 = \left[ \begin{array}{cc}
1 -i  & 1  \\
1  & 1 +i
\end{array}\right]
~~~~\mbox{and}~~~~
\bP_2 = \left[ \begin{array}{cc}
1  & 0  \\
-2i  & 1
\end{array}\right]
\eeq
has a limit set $\Lambda_{\Gamma}$ which is the complement
of the interiors of the circles in the Apollonian packing $\sP_{\sD_0}$
associated to the Descartes configuration $\sD_0$ above. 
\end{prop}

This limit set,  ``the glowing gasket'',  is illustrated in Figure 7.3 in 
Mumford, Series and Wright \cite[Chapter 7]{MSW02}. 

The Schottky group $\Gamma_{S}(\sD_0)$ is very special;
it admits no deformations other than conjugacy, and
under conjugacy $\fg \Gamma_{S}(\sD_0)\fg^{-1} =\Gamma_{S}(\sD) $
with $\sD = \fg(\sD_0)$.
 Recall that
a  matrix $\bM \in PSL(2, \CC)$ is
{\em parabolic} if $Trace(\bM) = \pm 2$; its 
action on the Riemann sphere has a single fixed point.
Here $\bP_1, \bP_2$ are parabolic, as is
\beql{1001b}
[\bP_1, \bP_2] :=  \bP_1 \bP_2 \bP_1^{-1} \bP_2^{-1} =  
\left[ \begin{array}{cc}
-1 -2i  & 2i  \\
-2i  & -1 +2i
\end{array}\right].
\eeq
The Schottky group $\Gamma_{S}(\sD_0)$ is
characterized up to conjugacy
by the property that it is generated by 
two elements $\bP_1, \bP_2$
which are parabolic,  such that their commutator 
$[\bP_1, \bP_2] = \bP_1 \bP_2 \bP_1^{-1} \bP_2^{-1}$ is
also parabolic (\cite[pp. 207--208]{MSW02}).

Now consider a the  group 
$G_{\sA}({\sD_0}) = \langle \fs_1, \fs_2, \fs_3, \fs_4 \rangle$
of M\"{o}bius transformations 
corresponding to the Apollonian group for
the Descartes configuration $\sD_0$, see \eqn{4311}.
Here $\bs_i$ is inversion in the i-th circle of the
dual Descartes configuration $\sD^{\perp}$, which is the 
circle passing through the three intersection points in $\sD$
not touching $C_i$.
The group $G_{\sA}({\sD_0})$ lifts to a group $\Gamma_{A}(\sD_0)$
of isometries of $\HH^3$. 
We describe $G_{\sA}({\sD_0})$  for  the case that
$\sD_0$ is the Descartes configuration in Figure~\ref{B1-special}
with circles numbered as in Figure~\ref{B1-special}(b) ; 
all other
$G_{\sA}(\sD)$ are related to it by conjugacy in
the M\"{o}bius group.
An inversion in a circle $C$ of radius $r$
centered at $\bz_0$ is the anti-holomorphic map 
$$
\fj_{C}(z)= z_0 + \frac{r^{2}}{\overline{z} - \overline{z_0}}
= \frac{ z_0 \bar{z} + r^{2} - |z_0|^2}{\bar{z} - \bar{z}_0}.
$$
We can write  the 
group generators  as $\fs_j =\fp_j \circ \fc$, 
where $\fp_j \in PSL(2, \cc)$  is holomorphic and $\fc$ denotes
complex conjugation, with
$$
\fp_1 = \left[ \begin{array}{cc}
1-i & -i \\
i & 1+i
\end{array}\right],~~~
\fp_2 = \left[ \begin{array}{cc}
1+i  & -i \\
i &  1-i
\end{array}\right],~~~
\fp_3 = \left[ \begin{array}{cc}
1 & 0 \\
4i  & 1
\end{array}\right],~~~
\fp_4 = \left[ \begin{array}{cc}
1 & 0 \\
0 & 1
\end{array}\right].
$$
Since $\fp_4 = \bI$,
the  subgroup of holomorphic elements of $G_{\sA}(\sD)$ is
\begin{eqnarray}
G_{\sA}^2(\sD) &= & \langle \fp_1\overline{\fp}_2, \fp_1\overline{\fp}_3, 
\fp_1 \overline{\fp}_4,
\fp_2 \overline{\fp_3}, \fp_2\overline{\fp_4}, \fp_3\overline{\fp_4},
 \fp_2\overline{\fp}_1, \fp_3\overline{\fp}_1, \fp_4 \overline{\fp}_1,
\fp_3 \overline{\fp_2}, \fp_4\overline {\fp}_2, \fp_4\overline{\fp_3}
 \rangle \nonumber \\
&=& \langle \fp_1, \fp_2, \fp_3 \rangle,
\end{eqnarray}
where we used the fact that $\overline{\fp_j} = \fp_j^{-1}$.

There are clearly relations
between the 
Schottky  group $\Gamma_{S}(\sD_0)$ and the 
 group $G_{\sA}^2(\sD_0)$ of holomorphic elements of $G_{\sA}(\sD)$
Indeed the group $G_{\sA}^2(\sD_0)$ is generated by three parabolic
elements $\fp_1, \fp_2, \fp_3$, and these are related to the three
parabolic elements $P_1, P_2, [P_1, P_2]$ 
in $\Gamma_{S}(\sD)$ by the relations
$$
\fp_3 = (P_2)^{-2},~~~~~\fp_1^{-1}\fp_2= -(P_1)^{-2},~~~~~
\fp_2^2 = -[P_1, P_2].
$$
However these groups are not exactly the same. It remains to determine
a more precise relation between them,  explaining how they are 
give rise to the  same Apollonian packing. Note that both these
groups lifted to discrete isometry groups of $\HH^3$
have fundamental domains in $\HH^3$ of infinite volume.

{\tt
\begin{tabular}{lllll}
email: & graham@ucsd.edu \\
& lagarias@umich.edu \\
  & colinm@research.avayalabs.com \\
 & allan@research.att.com \\
& cyan@math.tamu.edu
\end{tabular}
 }

\end{document}